\documentstyle[verbatim,11pt,amscd]{amsart}

\theoremstyle{plain}
\newtheorem{Thm}{Theorem}
\newtheorem{Lem}[Thm]{Lemma}
\newtheorem{Cor}[Thm]{Corollary}

\newtheorem{Prop}[Thm]{Proposition}
\newtheorem{Nonsep}{Theorem~\ref{nonsep}}

\newtheorem{Sep}{Theorem~\ref{sep}}

\newtheorem{Cornonumb}{Corollary}

\theoremstyle{definition}
\newtheorem{Def}{Definition}

\theoremstyle{remark}
\newtheorem{Rem}{Remark}

\newtheorem{Rems}{Remarks}

\def\IC{\Bbb C}

\def\cp{\hbox{${\Bbb C} P^2$}}
\def\cpb{\hbox{$\overline{{\Bbb C}P^2}$}}

\def\sx{\hbox{$S^2 \times S^2$}}
\def\snx{\hbox{$S^2 \widetilde{\times} S^2$}}

\def\a{\alpha}
\def\b{\beta}
\def\d{\delta}
\def\e{\epsilon}
\def\g{\gamma}
\def\S{\Sigma}
\def\i{\iota}
\def\s{\sigma}
\def\vp{\varphi}
\def\bdy{\partial}

\input epsf

\begin{document}
\title[]{Hyperelliptic Lefschetz fibrations and branched covering 
spaces}
\author{Terry Fuller}
\address{Department of Mathematics\\University of California\\
Irvine, CA 92697}
\curraddr{School of Mathematics\\
Institute for Advanced Study\\
Olden Lane\\
Princeton, NJ 08540}

\email{fuller@@math.ias.edu}
\subjclass{57R99, 57M12}
\thanks{The author was supported by NSF
grant DMS 97-29992.}

\begin{abstract}
Let $M$ be a smooth 4-manifold which admits a
relatively minimal hyperelliptic
genus $h$ Lefschetz fibration over $S^2$.
If all of the vanishing cycles for this fibration
are nonseparating curves, then we show that $M$ 
is a 2-fold cover of an $S^2$-bundle over $S^2$,
branched over an embedded surface.
If
the collection of vanishing cycles for this fibration includes
$\s$ separating curves,
we show that $M$ is the
relative minimalization of a 
Lefschetz fibration constructed as
a 2-fold branched cover of
$\cp \# (2\s+1) \cpb,$
branched over an embedded surface.

\end{abstract}

\maketitle
\setcounter{section}{-1}
\section{Introduction}

As a result of the crucial role currently played by symplectic
manifolds in smooth 4-manifold topology,  the notion of a
smooth Lefschetz fibration has recently taken on renewed
importance. It is now known that the existence of such a fibration
provides a purely topological description of symplectic 4-manifolds.
In particular, by a theorem of Donaldson \cite{d} it is known that
every symplectic 4-manifold admits a smooth Lefschetz pencil, which
can be blown up to yield a Lefschetz fibration over $S^2$.
Conversely, Gompf \cite{gs} has shown that every smooth 4-manifold
admitting a Lefschetz fibration (with a few understood exceptions)
is in fact symplectic, with
symplectic fibers.

A genus $h$ Lefschetz fibration (see below for a careful definition)
is a smooth fibration $M\rightarrow C$ of an oriented
4-manifold $M$ over an oriented 2-manifold $C$ by possibly
singular genus $h$ 2-manifolds, where every singular point
locally has the holomorphic description of two complex lines
meeting in a transverse double point. This strong condition
completely determines the local topology of the fibration,
with the relevant data phrased in terms of the vanishing
cycles of complex algebraic geometry.
Intuitively, a vanishing cycle is a simple closed curve on a nonsingular
fiber which is gradually shrunk to a point as one approaches the 
singular fiber. The diffeomorphism classification of regular 
neighborhoods of possible singular fibers in genus $h$
Lefschetz fibrations is particularly easy to state:
such neighborhoods are determined by whether their corresponding
vanishing cycles are nonseparating or separating curves, and in the
latter case by the genera of the separated surfaces.

A natural area of investigation is to study the extent to
which an arbitrary smooth Lefschetz fibration globally
resembles a holomorphic one. It is conjectured that any
hyperelliptic smooth Lefschetz fibration featuring only
nonseparating vanishing cycles is holomorphic \cite{st}.
It is known that all holomorphic hyperelliptic
Lefschetz fibrations can be obtained as 2-fold covers
of ruled surfaces, branched over smoothly embedded
curves \cite{p1}. Our main theorems demonstrate that this feature
extends to smooth Lefschetz fibrations over $S^2$.

{\Nonsep 
Let $M\to S^2$ be a relatively minimal hyperelliptic genus $h$
Lefschetz fibration, and assume that all of the
vanishing cycles of this fibration are nonseparating
curves.
Then $M$ is
a 2-fold cover of \sx or \snx,
branched over an embedded surface.}

{\Sep
Let $M\to S^2$ be a relatively minimal hyperelliptic genus $h$
Lefschetz fibration whose collection of vanishing cycles
includes separating curves.
Then $M\to S^2$ is 
the relative minimalization of a Lefschetz fibration
$M^{\prime}\to S^2$, and $M^{\prime}$ is a
2-fold cover of a rational surface,
branched over an embedded surface.
(The terms ``relatively minimal'' and 
``relative minimalization'' are defined below.)}

Since all genus 2 Lefschetz fibrations are hyperelliptic,
we immediately have the following corollaries.

{\Cornonumb
Let $M\to S^2$ be a relatively minimal genus 2
Lefschetz fibration, and assume that all of the
vanishing cycles of this fibration are nonseparating
curves.
Then $M$ is
a 2-fold cover of \sx or \snx,
branched over an embedded surface.}

{\Cornonumb
Let $M\to S^2$ be a relatively minimal genus 2
Lefschetz fibration whose collection of vanishing cycles
includes separating curves.
Then $M\to S^2$ is 
the relative minimalization of a Lefschetz fibration
$M^{\prime}\to S^2$, and $M^{\prime}$ is a
2-fold cover of a rational surface,
branched over an embedded surface.}

Our theorems have been proven independently by Siebert and
Tian, who demonstrate in addition that the branched coverings
of the theorems are branched over symplectically
embedded surfaces \cite{st}.
Our corollaries have also been proven
independently
by Smith \cite{s}.

Following preliminary sections on Lefschetz fibrations
and branched coverings, we prove Theorem~\ref{nonsep} in
Section 3 and Theorem~\ref{sep} in Section 4. 
The nature of the proof makes it possible to 
explicitly describe the branch sets obtained
in our theorems, a fact that we will illustrate by example
in Section 5. 

Of more interest, perhaps, than the theorems themselves
is the method of proof, which in each case is 
to fashion
the branched cover by hand, beginning with a
handlebody description of $M$ inherited from its
structure as a hyperelliptic Lefschetz fibration. Our construction
in Section 4 is motivated by work of Persson, who 
constructed genus 2 holomorphic Lefschetz fibrations as
double covers of ruled surfaces, and masterfully introduced
infinitely close triple point singularities into the branch
sets to study questions of complex surface geography \cite{p1}.
The proofs given here generalize Persson's work to the smooth
category, and to fibrations of any genus. We do this by considering
deformations and resolutions of infinitely close $n$-tuple point
singularities, which is discussed in Section 6.
We demonstrate that resolutions of infinitely close 
singularities yield Lefschetz fibrations with separating vanishing 
cycles, a phenomenon which is known by algebraic geometers
to hold for holomorphic genus 2 Lefschetz fibrations \cite{p2}.
Thus another facet of this work is to provide a purely topological proof
of this correspondence for smooth Lefschetz fibrations of any genus.
This correspondence has also been observed by Siebert and Tian 
for Lefschetz fibrations of any genus \cite{st},
and by Smith for genus 2 Lefschetz fibrations \cite{s}.

{\bf Acknowledgments.} 
In earlier preprint versions of this work,
the author claimed versions of Theorem~\ref{nonsep}
and Theorem~\ref{sep} describing all 
(i.e. non-hyperelliptic) Lefschetz fibrations
as simple 3-fold covers.
It is a pleasure to acknowledge correspondence with
Bob Gompf and Bernd Siebert which brought to light
a gap in those arguments. Indeed, recent work of
Siebert and Tian demonstrates that such a result is
impossible \cite{st}. The author would also like
to thank Ivan Smith for valuable conversations.

\section{Lefschetz Fibrations}

{\Def Let $M$ be a compact, oriented smooth 4-manifold, and
let $C$ be a compact, oriented 2-manifold. A proper smooth
map $f:M\to C$ is a {\em Lefschetz fibration}
if 
\begin{itemize}
\item[(1)] 
each critical value $x_1,\ldots,x_{\mu}$
of $f$ lies in 
${\mathrm interior}(C)$; and
\item[(2)] about each $f^{-1}(x_i)$ and $x_i$, there are complex 
coordinate charts agreeing in orientation with the orientations
of $M$ and $C$, respectively, such that locally $f$ can be expressed
as $f(z_1,z_2)=z_1^2+z_2^2.$
\end{itemize}}

It is a consequence of this definition that
$$f|_{f^{-1}(C-\{x_1,\ldots,x_{\mu})}:
f^{-1}(C-\{x_1,\ldots,x_{\mu}\})\to C-\{x_1,\ldots,x_{\mu}\}$$
is a smooth fiber bundle over $C-\{x_1,\ldots,x_{\mu}\}$
with fiber diffeomorphic to a 2-manifold $\S_h$, and so
we also refer to $f$ as a {\em genus $h$ Lefschetz fibration.}
Two genus $h$ Lefschetz fibrations $f:M\to C$ and
$f^{\prime}:M^{\prime}\to C^{\prime}$ are {\em equivalent}
if there are diffeomorphisms $\Phi:M\to M^{\prime}$ and
$\phi:C\to C^{\prime}$ such that 
$f^{\prime}\Phi=\phi f.$
In this paper, we will always have either $C=D^2$, or
$C=S^2$.

Given a Lefschetz fibration $M\to C$, one may blow up
$M$ to produce a Lefschetz fibration 
$M\# \cpb \to C$.
This will have the effect of replacing a point in
one fiber with an embedded 2-sphere of self-intersection number $-1$.
If no fiber of $M\to C$ contains an embedded sphere
of square $-1$, we say that 
$M\to C$ is {\em relatively minimal.}
If the Lefschetz fibration $M\to C$ is obtained
from a Lefschetz fibration $M^{\prime} \to C$
by blowing down all spheres of square $-1$ found in
fibers, we say that $M\to C$ is the 
{\em relative minimalization} of $M^{\prime} \to C$.

If $f:M\to S^2$ is a smooth genus $h$ Lefschetz fibration,
then we can use the Lefschetz fibration to produce a handlebody
description of $M$. Let $M_0=M-\nu(f^{-1}(x))$, where
$\nu(f^{-1}(x))\cong \S_h\times D^2$ denotes a regular neighborhood 
of a nonsingular fiber $f^{-1}(x)$.
Then $f|_{M_0}:M_0 \to D^2$ is a smooth Lefschetz fibration.
We select a regular value $x_0 \in {\mathrm interior}(D^2)$
of $f$, an identification $f^{-1}(x_0)\cong \S_h$,
and a collection of arcs $s_i$ in
${\mathrm interior}(D^2)$ with each $s_i$
connecting $x_0$ to $x_i$, and otherwise disjoint
from the other arcs. We also assume that the critical
values are indexed so that the arcs $s_1,\ldots,s_\mu$
appear in order as we travel counterclockwise in a 
small circle about $x_0$.
Let $V_0,\ldots,V_n$ denote a collection of small disjoint
open disks with $x_i \in V_i$ for each $i$.

To build our description of $M$, we observe first
that $f^{-1}(V_0)\cong \S_h\times D^2$,
with $\bdy V_0\cong \S_h\times S^1.$
Enlarging $V_0$ to include the critical value
$x_1$, it can be shown that
$f^{-1}(V_0\cup \nu(s_1) \cup V_1)$
is diffeomorphic to $\S_h\times D^2$
with a 2-handle $H_1^2$ attached along a circle
$\g_1$ contained in a fiber
$\S_h\times pt. \subset \S_h\times S^1.$
Moreover, condition (2) in the definition of a
Lefschetz fibration requires that $H_1^2$
is attached with a framing $-1$ relative
to the natural framing on $\g_1$
inherited from the product structure of $\bdy V_0$.
(See \cite{gs} and \cite{k} for proofs of these non-trivial
statements, and for more on the topology of Lefschetz fibrations.)
For intuition, one should picture the singular fiber
$f^{-1}(x_1)$ as being obtained by gradually shrinking 
the circle $\g_1$
to a point using the disk obtained from the core of $H_1^2$;
the circle $\g_1$ is traditionally dubbed a {\it vanishing cycle.}
In addition, $\bdy(\S_h\times D^2\cup H_1^2)$ is
diffeomorphic to a $\S_h$-bundle over $S^1$ whose monodromy
is given by $D_{\g_1}$, a righthanded Dehn twist about
$\g_1$. 

Continuing counterclockwise about $x_0$,
we add the remaining critical values to our description, 
yielding that
$$M_0\cong f^{-1}(
V_0\cup \bigcup_{i=1}^{\mu} \nu(s_i)
\cup \bigcup_{i=1}^{\mu} V_i)$$
is diffeomorphic to $\S_h \times D^2 \cup 
\bigcup_{i=1}^\mu H_i^2$, where each $H_i^2$
is a 2-handle attached along a vanishing cycle $\g_i$
in a $\S_h$-fiber in $\S_h\times S^1$ with relative framing $-1$.
Furthermore, 
$$\bdy M_0 \cong
\bdy(\S_h\times D^2 \cup \bigcup_{i=1}^\mu H_i^2)$$
is a $\S_h$-bundle over $S^1$ with monodromy given
by the product $D_{\g_1}\cdots D_{\g_\mu}$.
(This product should be read left to right, thus
$D_{\g_1}\cdots D_{\g_\mu}=
D_{\g_\mu}\circ\cdots\circ D_{\g_1}.$)
Since also $\bdy M_0 \cong \S_h\times S^1$, the 
global monodromy 
$D_{\g_1}\cdots D_{\g_\mu}$
is necessarily isotopic to the identity. 

Finally, to extend our description
of $M_0$ to $M$, we reattach a regular neighborhood 
$\nu(f^{-1}(x))$.
The total space $M$ can then be described as
$$M=\S_h\times D^2 \cup \bigcup_{i=1}^\mu H_i^2
\cup \S_h\times D^2,$$
where the final $\S_h\times D^2$ is attached
via a $\S_h$-fiber preserving map of the boundary.
The result is unique up to equivalence
for $h\geq 2$ \cite{k}.

Although the description of the monodromy
corresponding to each individual critical value $x_i$
as a Dehn twist depends on the choice of arc $s_i$,
other choices of arcs (and of the central identification
$f^{-1}(x_0)\cong \S_h$) do not change the Lefschetz
fibration on $M_0$, up to equivalence. 

We may also reverse this description to construct
Lefschetz fibrations, as follows.
Starting with a cyclically ordered collection of curves
$(\g_1, \ldots, \g_\mu)$ 
on $\S_h$, with
$D_{\g_1}\cdots D_{\g_\mu}$
isotopic to the identity, we may construct
a smooth 4-manifold $M_0$ 
by forming 
$\S_h\times D^2 \cup \cup_{i=1}^\mu H_i^2,$
where the 2-handles $H_i^2$ are attached in
clockwise order along
$\g_i \subset \S_h\times \{pt.\}
\subset \S_h\times S^1 = \bdy(\S_h\times D^2),$ 
each with framing $-1$ relative to the product framing
on $\S_h\times S^1$.
Since $\bdy M_0$ is a trivial
$\S_h$-bundle over $S^1$, we may then attach
$\S_h\times D^2$ to $M_0$ using a $\S_h$-fiber
preserving map of their boundaries, producing
a genus $h$ 
Lefschetz fibration over $S^2.$

{\bf Hyperelliptic Lefschetz fibrations.}
We say that a homeomorphism $\varphi:\S_h \to \S_h$
is {\it symmetric} if it commutes with the hyperelliptic
involution $\i: \S_h \to \S_h$.
We let ${\cal{H}}_h$ denote the {\it hyperelliptic mapping class
group,} namely the subgroup of ${\cal{M}}_h$ of mapping
classes that commute with the class of $\i.$

{\Def Let $M\to S^2$ be a genus $h$ Lefschetz fibration
whose global monodromy is given by the collection
$(\g_1, \ldots, \g_\mu)$ of vanishing cycles. Then
$M\to S^2$ is a {\it hyperelliptic Lefschetz fibration}
if there exists a $\varphi \in {\cal{M}}_h$ such that
$\varphi \circ D_{\g_i} \circ \varphi^{-1} \in {\cal{H}}_h$
for all $1\leq i \leq \mu.$}

The condition for a Lefschetz fibration to be hyperelliptic
is equivalent to the supposition that there is an identification
of a regular fiber $f^{-1}(x_0)$ with $\S_h$ so that
each of the vanishing cycles (up to isotopy) satisfies
$\i(\g_i)=\g_i.$ We refer to any curve $\g$ in $\S_h$ with
$\i(\g)=\g$ as {\it symmetric.}

\section{Branched Covers}

Let
$\pi:\S_h\to S^2$ denote the usual
2-fold cover branched over $2h+2$ points
obtained as the quotient of $\i$,
and let $B$ denote the branch set of $\pi$.
An easy Euler characteristic argument gives the following.

{\Lem \label{curvelift} 
Let
$\g \subset \S_h$ be an embedded symmetric closed curve, and let
$\d=\pi(\g)$.
Then either
\begin{itemize}
\item[(a)] $\d$ is an arc with endpoints on $B$, and otherwise
disjoint from $B$, with
$\pi |\g : \g\to\d$ a 2-fold branched cover; or
\item[(b)] $\d$ is a simple closed curve disjoint from $B$, with
$\pi |\g : \g\to\d$ a 2-fold (unbranched) cover.
\end{itemize}
}

We note that case (a) in the above Lemma occurs when $\g$ is
a nonseparating curve on $\S_h$, and case (b) occurs when $\g$ is
separating.

{\Lem \label{homeolift}
Let $\g_1$ and $\g_2$ be embedded symmetric closed curves
on $\S_h$, and assume that 
both separate $\S_h$ into surfaces of the same genera.
Then there exists a symmetric homeomorphism
$\tilde{f}:\S_h\to \S_h$ and
a homeomorphism
$f:S^2\to S^2$ such that 
$\tilde{f}(\g_1)=\g_2$ and
the diagram
$$\begin{CD}
\S_h @>\tilde{f}>> \S_h\\
@V{\pi}VV          @VV{\pi}V\\
S^2  @>f>> S^2
\end{CD}$$
commutes.}

\begin{pf}
From Lemma~\ref{curvelift}, case (b), 
$\pi(\g_1)=\d_1$ and $\pi(\g_2)=\d_2$
are each simple closed curves in $S^2$
which, by assumption, separate $S^2$
into regions containing the same number
of points of $B$. Thus, we may find a 
homeomorphism $f:S^2\to S^2$ fixing the
set $B$ with $f(\d_1)=\d_2.$
The mapping class group
${\cal M}_{0,2h+2}$ is known to be
generated by ``disk twists'' about embedded
arcs with endpoints on $B$, that is, by
homeomorphisms which rotate a small disk neighborhood
of the arc by $180^{\circ},$ exchanging the endpoints \cite{b}. 
Thus we may (after an isotopy) factor $f$ as the product
of such disk twists. This product realizing $f$ on $S^2$ 
will lift via $\pi$ to a product $\tilde{f}$ of
Dehn twists about symmetric curves, hence to
a symmetric homeomorphism with 
$\tilde{f}(\g_1)=\g_2$.
\end{pf}

\begin{figure}[h]
$$ \epsfbox{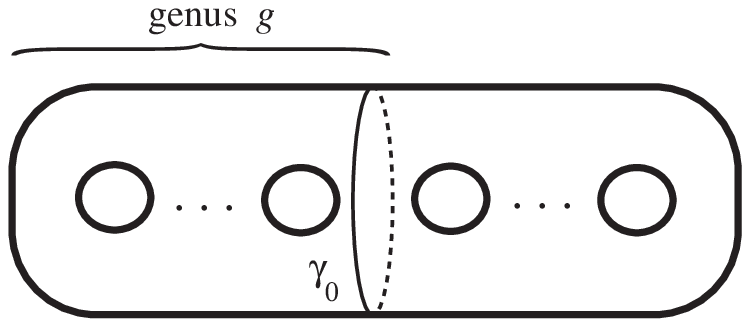} $$
\caption{} 
\label{gamma0}
\end{figure}

\begin{figure}[h]
$$\epsfbox{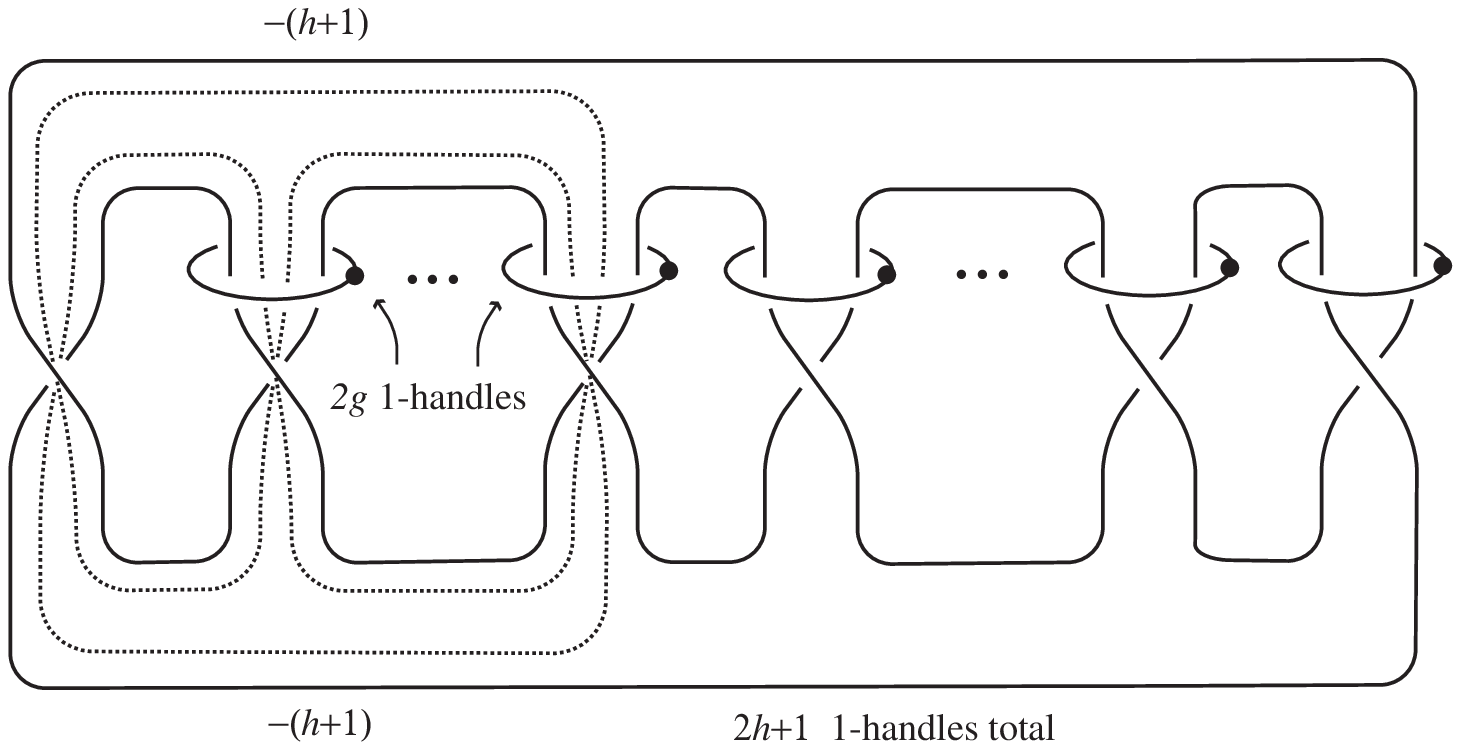}$$
\caption{} 
\label{gamma0lift}
\end{figure}

{\bf A separating curve}
Let $\g_0$ be the separating curve on $\S_h$ shown in Figure~\ref{gamma0},
and assume that $\g_0$ bounds a genus $g$ surface.
{\Lem \label{g0}
The 4-manifold described by the framed link in
Figure~\ref{gamma0lift} is diffeomorphic to
$\S_h\times D^2$. This diffeomorphism
maps the dashed curve in Figure~\ref{gamma0lift}
to $\g_0 \subset \S_h\times \{pt.\}
\subset \S_h\times S^1 =
\bdy(\S_h\times D^2).$
Furthermore, the product framing on this curve
in $\S_h\times S^1$ agrees with the $0$-framing
on this curve (defined by a Seifert surface in $S^3$.)
}

\begin{pf}
We will prove the lemma by drawing $\S_h\times D^2$
as the total space of the 2-fold branched cover
$\pi\times id: \S_h\times D^2 \to S^2\times D^2,$
using the algorithm of \cite{ak}.
Figure~\ref{h1} shows the base $S^2\times D^2$
of this cover, with the branch set visible as the
$2h+2$ disjoint horizontal disks $B\times D^2$.
\begin{figure}[h]
$$\epsfbox{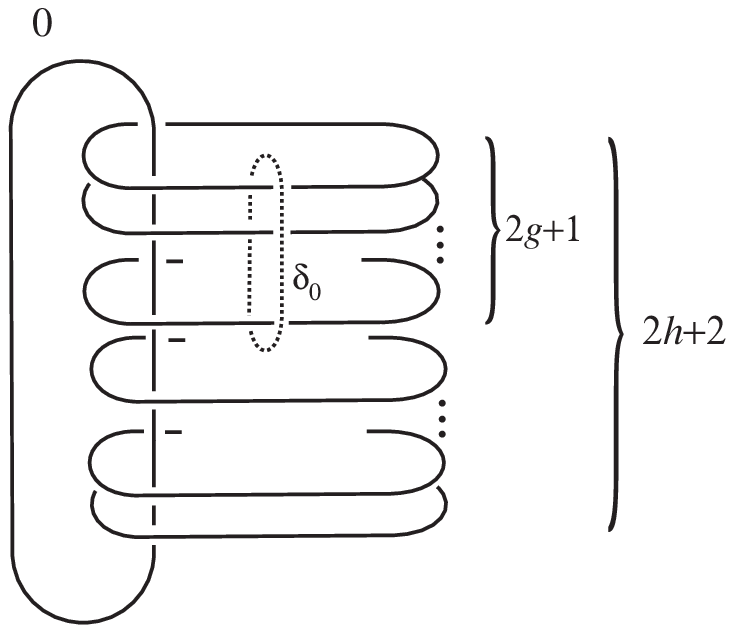}$$
\caption{} 
\label{h1}
\end{figure}
In this Figure, $\pi\times id$ restricted to the boundary
is easily visualized: picturing the standard
description of the complement of the $0$-framed
unknot as $D^2\times S^1,$ with the $S^1$ factor
given by a meridianal circle, we see the $2h+2$
branch points in each $D^2\times \{pt.\}$ 
on the boundary of each meridianal disk.
The curve $\g_0$ is the lift
under $\pi$ of the simple closed curve $\d_0$
shown in Figure~\ref{2fold}.
\begin{figure}[h]
$$\epsfbox{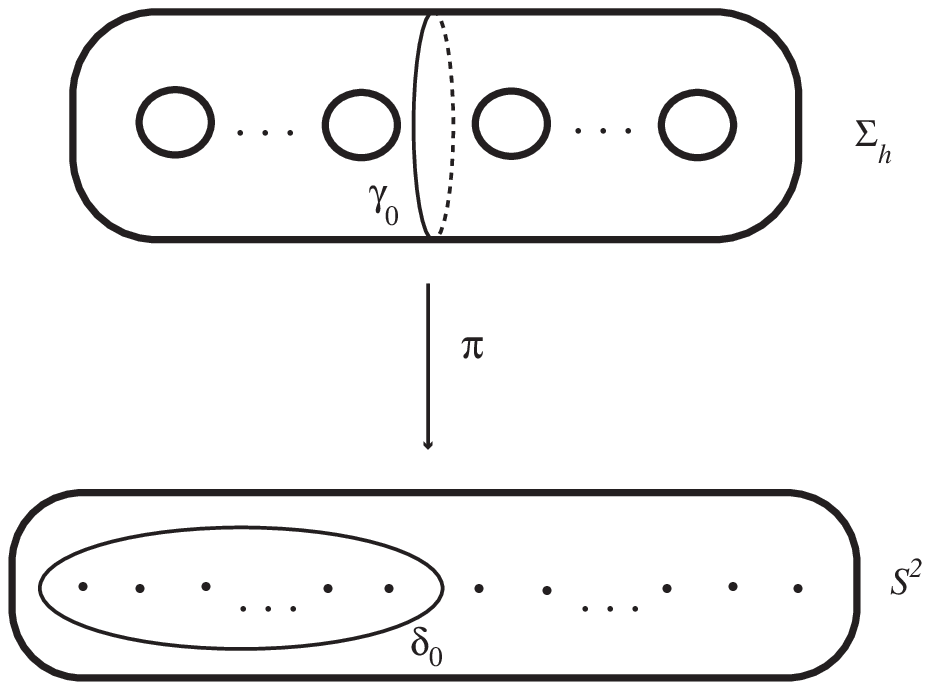}$$
\caption{} 
\label{2fold}
\end{figure}
The dashed curve in Figure~\ref{h1}
represents
$\d_0 \subset S^2\times \{pt.\} \subset S^2\times S^1=
\bdy(S^2\times D^2).$

To form the 2-fold branched cover shown in Figure~\ref{gamma0lift},
we imagine beginning with the unbranched 2-fold cover given by 
two disjoint copies of Figure~\ref{h1}.
The procedure of \cite{ak} is to push the interior of the branch
set into the interior of $S^2\times D^2$, cut along the track of 
this isotopy, and attach handles to glue the
copies together. More specifically, gluing the 0-handles of the
two copies together along the $2h+2$ disks of the branch set
in Figure~\ref{h1} results in the attachment of $2h+1$
1-handles, i.e. the row of dotted circles in Figure~\ref{gamma0lift}.

We must also lift the 0-framed 2-handle in Figure~\ref{h1}
to the 2-fold cover. Since its attaching circle intersects the branch
set, it will lift to two 2-handles obtained by cutting the attaching
circles in each copy of Figure~\ref{h1} where they intersect
the branch set, and joining the endpoints together as we perform
the gluings in the previous paragraph. Hence it lifts to the
two 2-handles shown in Figure~\ref{gamma0lift}.
We may calculate the framings of these lifts as in \cite{ak}:
letting $\lambda_1$ and
$\lambda_2$ denote the 2-dimensional homology classes
defined by these lifts of the 
2-handle in
Figure~\ref{h1}, we must have
$(\lambda_1 + \lambda_2)^2 = 2(0)=0.$
Expanding, and making use of symmetry, we get
$2\lambda_1^2 + 2\lambda_1 \cdot \lambda_2 =0$.
Reading $\lambda_1 \cdot \lambda_2=h+1$ directly
from Figure~\ref{gamma0lift}, we can then solve 
$\lambda_1^2=-(h+1)=\lambda_2^2.$

Finally, to see that the dashed curve in Figure~\ref{gamma0lift}
is $\g_0$, we simply note that we can draw the lift of $\d_0$
exactly as we lifted the the attaching circles of the 2-handle,
with the outcome as pictured.
To visualize the product framing on the dashed curve in
Figure~\ref{gamma0lift}, we can lift a curve parallel
to $\d_0$ in a different 
$S^2\times \{pt.\}$ fiber.
This parallel curve will lift to a curve which lies on
the obvious Seifert surface for the dashed curve
in Figure~\ref{gamma0lift}.
This confirms that the product framing and the
$0$-framing agree.
\end{pf}

{\section{Nonseparating Vanishing Cycles}

We are now ready for our first theorem.

{\Thm \label{nonsep}
Let $M \to S^2$ be a hyperelliptic genus $h$ 
Lefschetz fibration and assume that all
of the vanishing cycles of this fibration are nonseparating
curves. Then $M$ is a 2-fold cover of \sx or \snx,
branched over an embedded surface.
}

{\Rems It has long been known that every elliptic ($h=1$)
Lefschetz fibration over $S^2$ 
with at least one singular fiber
may be obtained as a 2-fold branched cover
of \sx. This follows from Moishezon's classification of
elliptic Lefschetz fibrations \cite{ms}.

It is known that a genus 2 Lefschetz fibration
where all vanishing cycles are about nonseparating
curves must have $10n$ singular fibers, for some 
integer $n\geq1$ \cite{ma}. Smith has shown
that every genus 2 Lefschetz fibration
is a 2-fold branched cover of \sx when $n$
is even, and is a 2-fold branched cover of 
\snx when $n$
is odd \cite{s}.}

\begin{pf}
From Section 1, we can describe $M$ as a handlebody
$$M=\S_h\times D^2 \cup \bigcup_{i=1}^\mu H_i^2
\cup \S_h\times D^2,$$
where each 2-handle $H_i^2$ is
attached to a symmetric vanishing cycle $\g_i$ with framing $-1$ 
relative to the induced 
framing coming from $\S_h \times S^1.$ 
Let $M_0$ denote the submanifold
$M_0=\S_h\times D^2\cup \bigcup_{i=1}^\mu H_i^2
\subset M$.
We begin with the 2-fold cover
$\pi\times id: \S_h\times D^2 \to S^2 \times D^2$
whose branch
set the $2h+2$  
disks $B\times D^2$ in Figure~\ref{h1}. 
Our first task is to extend this branched covering to $M_0$.

Since each vanishing cycle $\g_i$ is a 
nonseparating symmetric curve in
$\S_h \times \{pt.\} \subset \S_h \times S^1$,
we can apply case (a) of Lemma~\ref{curvelift} to each one to produce
an arc $\d_i \subset S^2\times \{pt.\}$. For each $\g_i$,
we then add a band to the branch set of Figure~\ref{h1}
whose core is $\d_i$, and which differs from the
band $\d_i\times (pt.-\epsilon,pt.+\epsilon) \subset
S^2 \times S^1$ by one lefthanded half twist. See Figure~\ref{base2}.
\begin{figure}
\leavevmode
$$\epsfbox{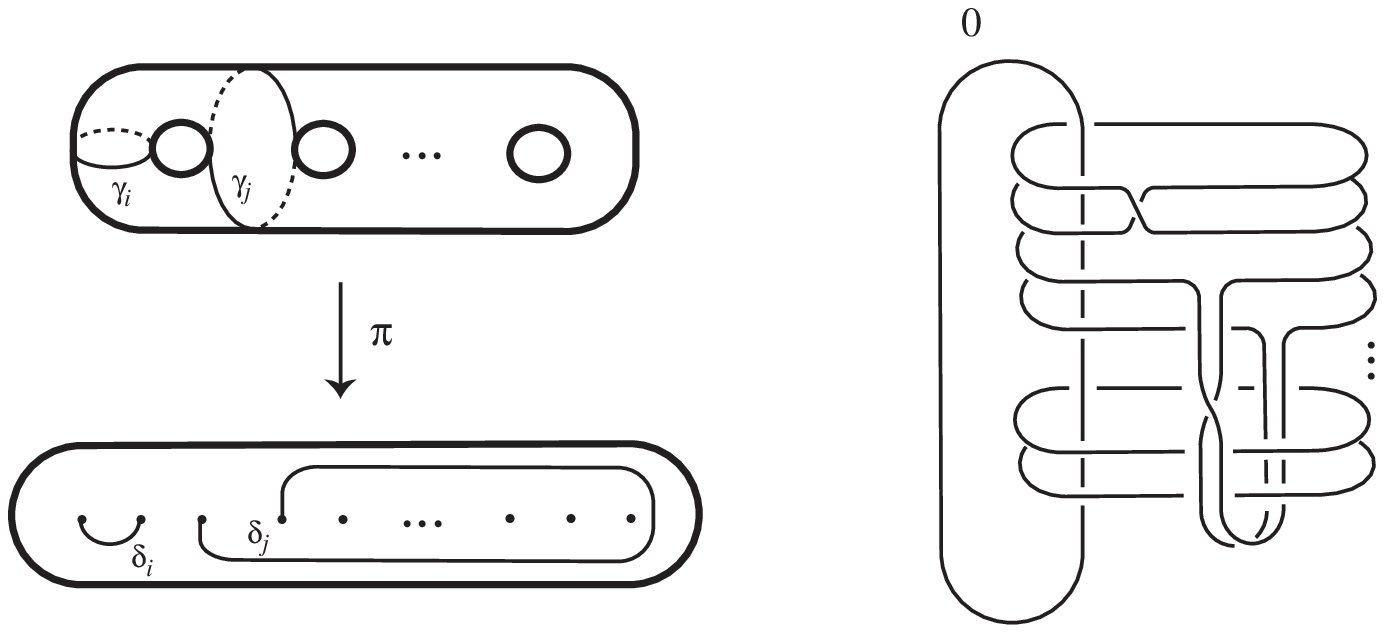}$$
\caption{} \label{base2}
\end{figure}
This new surface in $S^2 \times S^1$
may contain ribbon singularities,
as shown, and following \cite{mo} we refer to this surface as
a {\em ribbon manifold}. 
We may push the interior 
of this ribbon manifold into the interior of $\S^2 \times
D^2$ to produce an embedding.

{\Prop \label{ext}
The 2-fold cover of $S^2 \times D^2$
branched over the ribbon manifold constructed above
(with its interior pushed into the interior of
$S^2 \times D^2$) is diffeomorphic to $M_0$.}

\begin{pf*}{Proof of Proposition~\ref{ext}}
This follows immediately from the proof of Theorem~6 of
\cite{mo}, modified slightly for our setting.
For completeness, we summarize the argument and refer
the reader to \cite{mo} for additional details.
From Lemma~\ref{curvelift}, a regular neighborhood 
$\nu(\d_i)$ in $S^2\times S^1$ lifts under $\pi\times id$ to the 
regular neighborhood $\nu(\g_i)$ in $\S_h\times S^1$.
We attach 2-handles $H_i^2$ to $\S_h\times D^2$
with attaching circles $\g_i$ and relative framing $-1$ via
$h_i:(\bdy D^2\times D^2)_i \to \nu(\tilde{\g_i})$.
Additionally, we define $g_i:H_i^2 \to H_i^2/V$
to be the 2-fold cover induced from the involution
$V:D^2\times D^2 \to D^2\times D^2$ which is reflection in
$D^1 \times D^1$. 
We may assume, after an isotopy, that the involution
$h_i \circ V \circ h_i^{-1}$ on $\nu(\g_i)$ agrees with
the involution $\pi\times id$ on $\nu(\g_i)$,
hence this involution
extends over $H_i^2$.
We can then form the map
$$\pi\times id \cup {\textstyle \bigcup}_i g_i$$
from $$\S_g\times D^2\cup {\textstyle \bigcup}_{h_i}H_i^2$$ 
to $$S^2\times D^2 \cup  
{\textstyle \bigcup}_{(\pi\times id) \circ h_i \circ g_i^{-1}} H_i^2/V.$$
This map is a 2-fold cover
whose branch set is isotopic to the ribbon manifold above;
the half twist in the bands above are required due to the
assumption that each $H_i^2$ is attached with relative
framing $-1$.
However, the addition of the 4-balls $H_i^2/V$
to domain and range 
does not change the manifolds, 
so we have constructed
a 2-fold cover
$\S_h\times D^2\cup \bigcup_{h_i}H_i^2 \to
S^2\times D^2.$
\end{pf*}

Finally, it remains to extend the branched covering
giving $M_0$ to all of $M$. From Proposition~\ref{ext},
we have constructed a 2-fold cover
$\varphi:M_0\to S^2\times D^2$. 
By construction, $\bdy M_0$ is a $\S_h$-bundle
over $S^1$ with monodromy
$H=D_{\g_1}\cdots D_{\g_\mu}$, which is
isotopic to the identity.
Restricted to the boundary,
the image of $\varphi$ is an $S^2$-bundle over $S^1$
whose monodromy is induced from the word
$D_{\g_1}\cdots D_{\g_\mu}.$
Furthermore, $\vp$ restricted to each fiber of
$\bdy M_0$ is the 2-fold covering $\pi.$
(This may be seen explicitly.
The branched covering $\vp$ projects each Dehn twist
$D_{\g_i}$ to a disk twist about the arc $\d_i$.
The boundary of the branch set of $\varphi$
appears as the link $L$ in $S^2\times S^1$
that is the boundary of the ribbon
manifold branch surface in $S^2\times D^2$.
The half twist placed in each band above ensures that
$L$ is isotopic to a closed braid.
While traversing the $S^1$ factor of 
$S^2\times S^1$,
this braid records the motion of the
$2h+2$ branch points of $\pi$ induced from
the monodromy
$D_{\g_1}\cdots D_{\g_\mu}.$)

We have a $\S_h$-bundle equivalence
$\Psi:\bdy M_0 \to \S_h\times S^1$
obtained by expressing $\bdy M_0$ as 
$\S_h\times I/(x,1)\sim (H(x),0)$,
and using the isotopy between $H$ and the
identity to adjust the gluing.
By \cite{bh}, we may assume that this isotopy
is an isotopy through
homeomorphisms that preserve the fibers of $\pi$.
Thus, if $H_t$ denotes this
isotopy, then $\pi H_t$ projects to an
isotopy between $\pi H$ and the identity on $S^2$.
This gives a bundle equivalence
$\psi:S^2\times S^1 \to S^2\times S^1$
for which the diagram
$$\begin{CD}
\bdy M_0  
@>\Psi>> \S_h\times S^1\\
@V{\varphi}VV        @VV{\pi\times id}V\\
S^2\times S^1 @>\psi>> S^2\times S^1
\end{CD}$$
commutes.
Therefore,
$$\varphi\cup (\pi\times id): M_0 \cup_{\Psi^{-1}} (\S_g\times D^2)
\to (S^2\times D^2) \cup_{\psi^{-1}} (S^2\times D^2)$$
describes $M$ as a 2-fold cover of 
an $S^2$-bundle over $S^2$,
with branch set the closed surface
obtained by attaching the disks
$B\times D^2 \subset S^2\times D^2$ to the ribbon 
manifold produced in Proposition~\ref{ext}.
\end{pf}

\section{Separating Vanishing Cycles}

{\bf A Model Fibration.}
The construction of the branched covers in the previous
section relied crucially on the assumption that all 
vanishing cycles were nonseparating.
In this section, we wish to extend that constuction
to achieve arbitrary hyperelliptic Lefschetz fibrations, which may
include separating vanishing cycles, as 2-fold branched
covers. To do that, it will be necessary to have a local
model of how a Lefschetz fibration with a separating vanishing
cycle can arise from a branched covering construction.

With this in mind, 
let $N_0\to D^2$ denote the genus $h$ Lefschetz fibration
with one singular fiber whose vanishing cycle is the curve
$\g_0$ pictured in Figure~\ref{gamma0}.
We again
consider
the 2-fold cover
$\pi\times id:\S_h\times D^2 \to S^2\times D^2,$
whose base is
shown in Figure~\ref{h1}.
We next would like to extend $\pi\times id$ to be
a branched cover of $S^2\times D^2\# 2\cpb$.
We attach a 2-handle along
$\d_0 \subset S^2\times \{pt.\}$ with framing $-1$, and we
attach another 2-handle along a meridian $\e_0$ to $\d_0$ with
framing $-2$. See Figure~\ref{b2}.
\begin{figure}[h]
$$\epsfbox{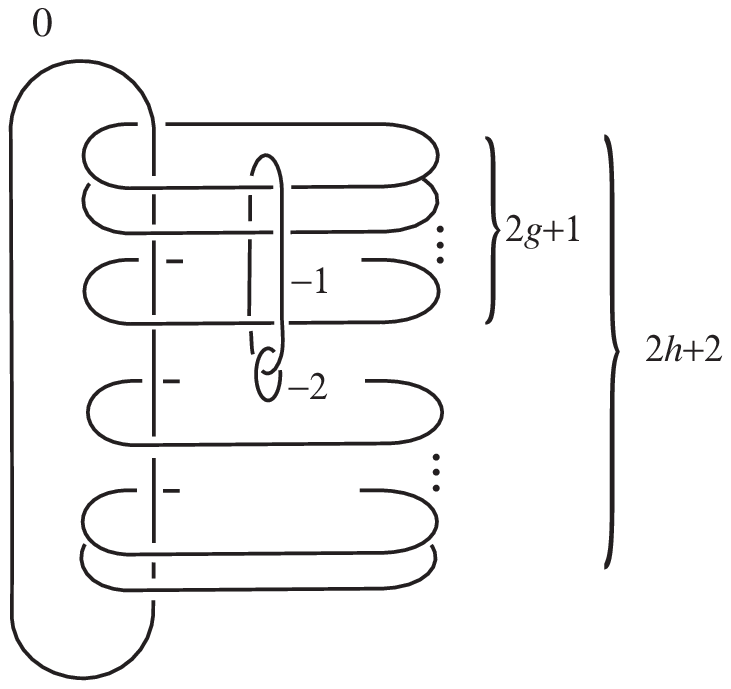}$$
\caption{} \label{b2}
\end{figure}
Blowing down the $-1$-framed 2-handle shows that
this manifold is 
diffeomorphic to $S^2\times D^2\# 2\cpb.$
We note for future reference that Figure~\ref{b2}
also reveals a boundary diffeomorphism
$\bdy(S^2\times D^2\# 2\cpb) \to S^2\times S^1$
obtained by blowing down twice. This boundary
diffeomorphism will put two full right-handed twists
in the upper strands of the boundary of the branch set.
Let $S$ denote the 2-sphere in $S^2\times D^2\# 2\cpb$
with $[S]^2=-2$ given by the $-2$-framed 2-handle in Figure~\ref{b2}.
To extend $\pi\times id,$ we define
$N_0^{\prime}$ to be the 2-fold cover of 
$S^2\times D^2\# 2\cpb$ ,
branched over
$B^{\prime}=B \cup S.$ 

{\Prop \label{model}
There is a non-relatively minimal
Lefschetz fibration $N_0^{\prime}\to D^2,$
and its relative minimalization is 
$N_0\to D^2$. In fact,
$N_0^{\prime}\cong N_0\# \cpb$.}

\begin{pf}
Figure~\ref{cover1} shows a framed link description of 
$N_0^{\prime}$,
drawn using the algorithm of
\cite{ak}. 
\begin{figure}[h]
\epsfxsize=6in 
$$\epsfbox{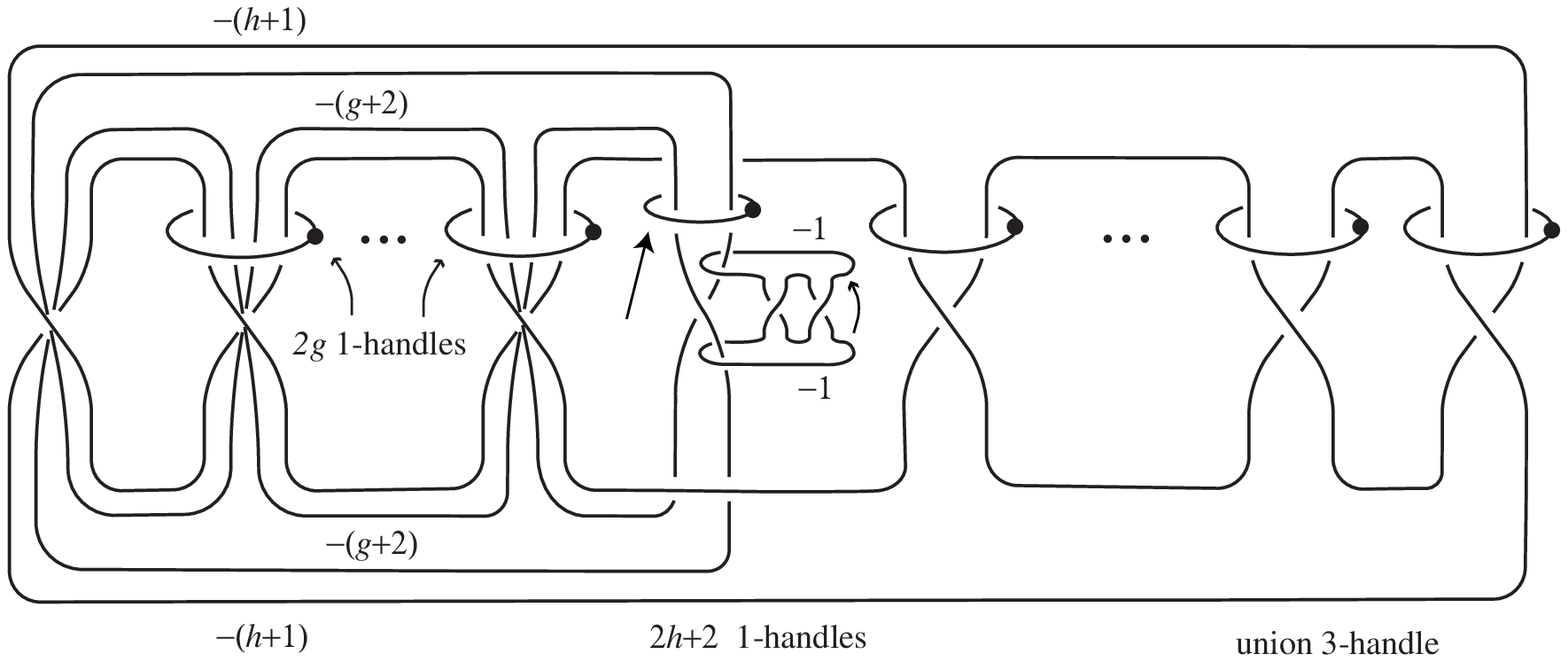}$$
\caption{} 
\label{cover1}
\end{figure}
The branch set $B^{\prime}$ is visible in Figure~\ref{b2} as
$2h+2$ long horizontal disks union
an embedded sphere $S$ of square $-2$ given
by the $-2$-framed 2-handle.
Hence $B^{\prime}$ has a handle decomposition
given by $2h+3$ 0-handles union a 2-handle
(coming from the core of the $-2$-framed
2-handle). 

To form the 2-fold branched cover shown in Figure~\ref{cover1},
we imagine beginning with the unbranched 2-fold cover
given by two disjoint copies of Figure~\ref{b2},
and we perform identifications as described in our
proof of Lemma~\ref{g0}.
Thus, gluing the 0-handles of
the two copies together along the 
$2h+2$ disks and the 0-handle of $S$ in Figure~\ref{b2}
results in the 
row of dotted circles in Figure~\ref{cover1}.
(The dotted circle marked with an arrow comes from $S$.)
As in Lemma~\ref{g0}, the 0-framed 2-handle in Figure~\ref{b2}
lifts to two 2-handles whose framings are 
$-(h+1).$ 
The $-1$-framed 2-handle in Figure~\ref{b2}
lifts to the two 2-handles of framing
$-(g+2)$, and the $-2$-framed 2-handle lifts
to the $-1$-framed 2-handles which link them.
Finally, we must glue the two copies of 
Figure~\ref{b2} along the 2-handle of $S$.
This results in the addition of a 3-handle to
Figure~\ref{cover1}, and completes the development of
that picture of $N_0^{\prime}.$

We next manipulate Figure~\ref{cover1} to see that
it has the properties claimed. We cancel the
1-handle marked with an arrow with one of the
$-(g+2)$-framed 2-handles. Sliding one of the
$-1$-framed 2-handles over the other
as indicated splits off a $0$-framed unknot,
which may be moved to be disjoint from the rest
of the picture, and used to cancel the 3-handle.
The result is Figure~\ref{cover3}.
\begin{figure}[h]
\epsfxsize=6in 
$$\epsfbox{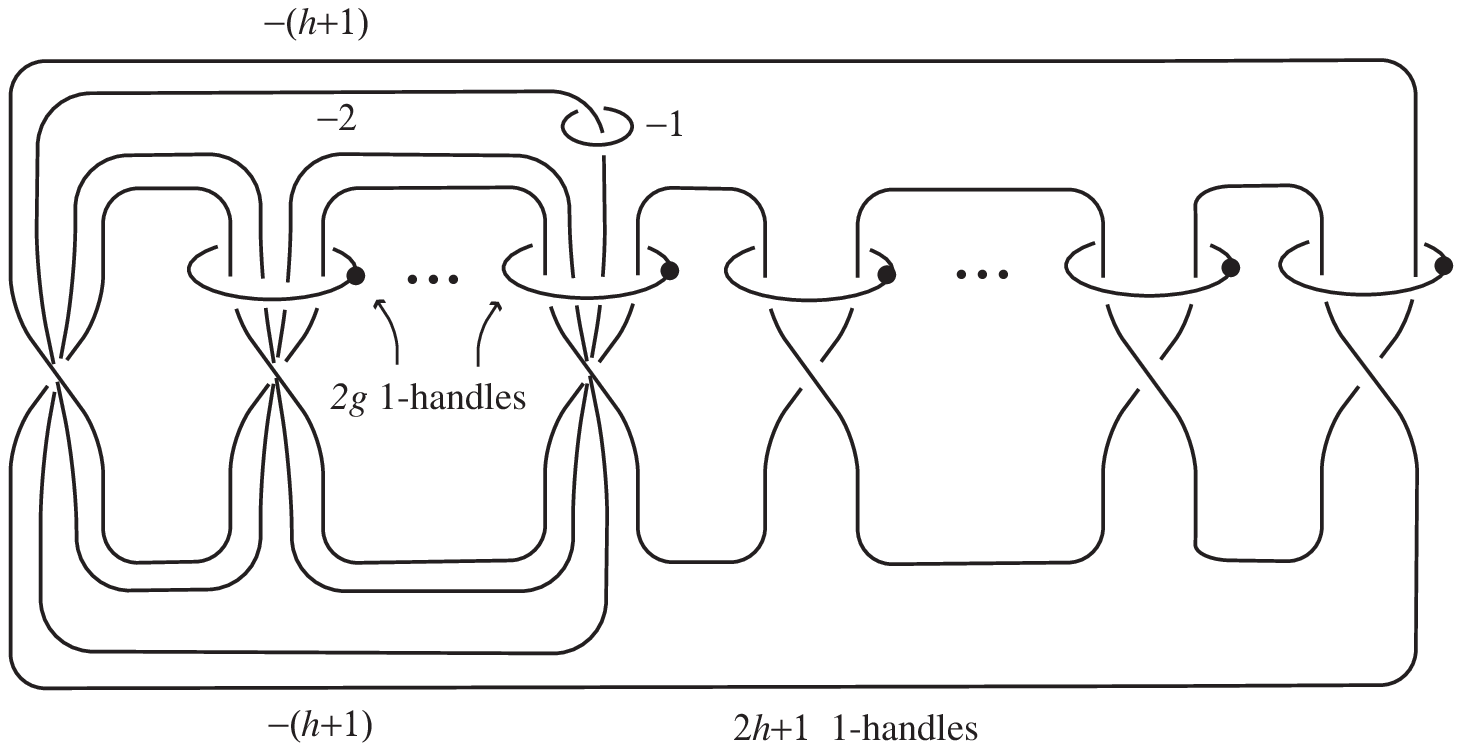}$$
\caption{} 
\label{cover3}
\end{figure}
Blowing down the $-1$-framed meridian to the
$-2$-framed 2-handle gives Figure~\ref{cover4}.
\begin{figure}[h]
\epsfxsize=6in 
$$\epsfbox{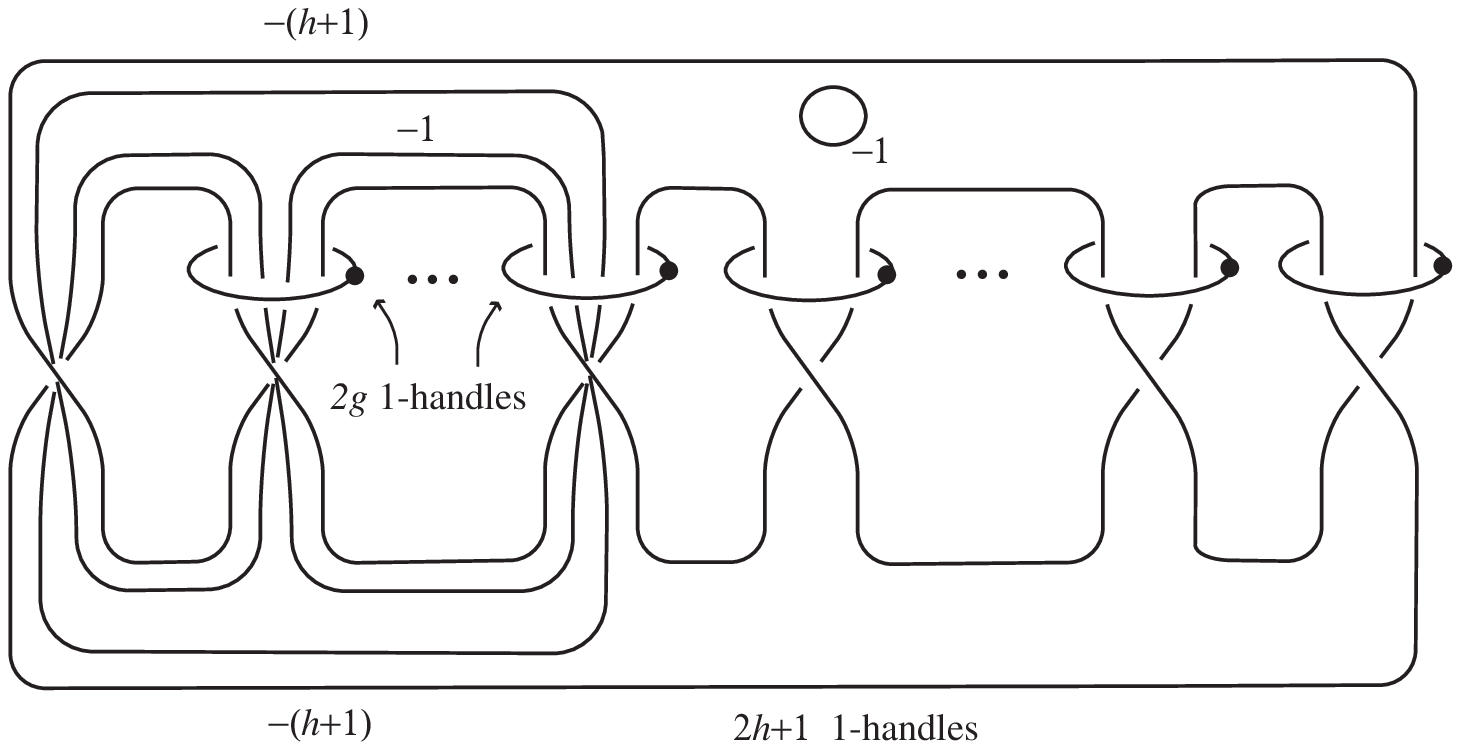}$$
\caption{} 
\label{cover4}
\end{figure}
Applying Lemma~\ref{g0}, we see that 
Figure~\ref{cover4} depicts 
$\S_h\times D^2$ with 2-handles
attached to curves 
$\g_0$ and $\a$
(where $\a$ bounds a disk in
$\S_h \times \{pt.\}$), each with relative framing $-1$.
This is
precisely the
Lefschetz fibration $N_0 \to D^2,$
blown up at three points.
\end{pf}

Examining Figure~\ref{cover4}, we recognize the
Lefschetz fibration 
$N_0\# \cpb \to D^2$ as having two vanishing
cycles $\g_0$ and $\a,$
where $\a$ bounds a disk in
$\S_h \times \{pt.\}.$
Analyzing the proof of Proposition~\ref{model}
more carefully, we have the following.

{\Cor \label{vccor} Let 
$\vp_0:N_0\# \cpb \to S^2\times D^2\# 2\cpb$ 
denote the
2-fold cover 
described in the proof of Proposition~\ref{model}.
Then $\vp_0^{-1}(\d_0)=\g_0$ and
$\vp_0^{-1}(\e_0)=\a.$}

\begin{pf}
We can record
$\vp_0^{-1}(\d_0)$ and $\vp_0^{-1}(\e_0)$ 
by reenacting the proof of Proposition~\ref{model},
recalling that $\d_0$ and $\e_0$ are the attaching 
circles of the $-1$- and $-2$-framed 2-handles,
respectively, in Figure~\ref{b2}.
Tracking the lifts of these handles through the
moves in the above proof gives the corollary.
\end{pf}

{\bf Separating Vanishing Cycles.}
We now state and prove our second main theorem.

{\Thm \label{sep}
Let $M\to S^2$ be a relatively minimal hyperelliptic
genus $h$ Lefschetz fibration
whose monodromy includes $\s$ separating vanishing cycles. Then
$M\to S^2$ is the relative minimalization of a Lefschetz fibration
$M\# \s\cpb\to S^2$, and $M\# \s\cpb$ is a 2-fold cover of
$\cp\# (2\s+1)\cpb$, branched over an
embedded surface.
}

\begin{pf}
From Section 1, we can represent $M$ as a handlebody
$$M=\S_h\times D^2 \cup \bigcup_{i=1}^\mu H_i^2
\cup \S_h\times D^2,$$
where each 2-handle $H_i^2$ is attached along a symmetric vanishing
cycle $\g_i$ with framing $-1$ relative to the induced
product framing on $\S_h\times S^1.$
Let $M_0=\S_h\times D^2\cup \bigcup_{i=1}^\mu H_i^2
\subset M$.
We again begin by forming
the 2-fold cover
$\pi\times id:\S_h\times D^2 \to S^2\times D^2.$

We would like to extend $\pi\times id$ to be
a branched covering with
$M_0 \# \s\cpb$
as its total space.
For the 2-handles $H_i^2$ whose vanishing cycles
are nonseparating curves, we may do this as in Section 3.
Thus, let $H_i^2$ be a 2-handle attached along a separating
vanishing cycle 
$\g_i \subset \S_h\times \{pt.\} \subset \S_h\times S^1,$
where $\g_i$ bounds a genus $g$ surface in $\S_h\times \{pt.\}.$
Let $N_i=\S_h\times D^2\cup H_i^2$, 
so that there is a
Lefschetz fibration $N_i\to D^2$ with 
one singular fiber whose vanishing cycle is $\g_i$.

From Lemma~\ref{homeolift}, we may find a symmetric
diffeomorphism 
$\tilde{f}:\S_h\to \S_h$ with 
$\tilde{f}(\g_0)=\g_i,$ 
and a diffeomorphism
$f:S^2\to S^2$ such that
the following diagram is commutative,
\begin{equation} \tag{*}
\begin{CD}
(\S_h\times D^2, \g_0) @>\tilde{f}\times id>> (\S_h\times D^2,\g_i)\\
@V{\pi\times id}VV          @VV{\pi\times id}V\\
(S^2\times D^2,\d_0)  @>f\times id>> (S^2\times D^2,\d_i)
\end{CD}
\end{equation}
where $\d_i=f(\d_0).$

The diffeomorphism 
$f\times id:S^2\times D^2 \to S^2\times D^2$
along the lower row extends to a diffeomorphism
$F:S^2\times D^2\# 2\cpb \to S^2\times D^2\# 2\cpb$,
where the range is Figure~\ref{h1} with a $-1$-framed
2-handle attached along $\d_i \subset S^2\times \{pt.\}$
and a $-2$-framed 2-handle attached along a meridian
to $\d_i$. The diffeomorphism of the upper row
extends to a diffeomorphism $\widetilde{F}: N_0 \to N_i$
(since it takes the vanishing cycle $\g_0$ to
$\g_i$), and further extends (by sending the
vanishing cycle $\a$
to a vanishing cycle which bounds a disk) to
a diffeomorphism
$\widetilde{F}: N_0\# \cpb \to N_i\# \cpb$.
Finally, we define $\vp_i=F\circ \vp_0 \circ 
\widetilde{F}^{-1}:N_i\# \cpb \to S^2\times D^2\# 2\cpb,$
which is a 2-fold branched cover.
This gives a commutative diagram
$$\begin{CD}
N_0\# \cpb @>\widetilde{F}>> N_i\# \cpb\\
@V{\vp_0}VV          @VV{\vp_i}V\\
S^2\times D^2\# 2\cpb  @>F>> S^2\times D^2\# 2\cpb.
\end{CD}$$
Furthermore, from the commutativity of
(*) and Corollary~\ref{vccor}, it follows that 
the 2-fold cover 
$\vp_i:N_i\# \cpb\to S^2\times D^2\# 2\cpb$
may be explicitly described by attaching a
$-1$-framed 2-handle to Figure~\ref{h1} along
$\d_i$, together with a $-2$-framed 2-handle along
a meridian to $\d_i$, and forming the 2-fold cover of
$S^2\times D^2\# 2\cpb$ branched over $B\cup S.$
Doing this for all $\s$ of the separating vanishing cycles
gives a 2-fold branched covering
$$\vp: M_0 \# \s\cpb
\to S^2\times D^2\# 2\s\cpb.$$

Finally, we must extend $\vp$ to give a branched covering
with total space $M\# \s\cpb$. 
We can modify the argument from Section 3 to the current
situation. Restricting $\vp$ to the boundary, 
we have a 2-fold cover
$$\bdy(M_0 \# \s\cpb)\to \bdy(S^2\times D^2\# 2\s\cpb)
\to S^2\times S^1,$$
where the latter map is the boundary diffeomorphism given by blowing
down. 
Let
$\omega:\bdy(M_0 \# \s\cpb)\to S^2\times S^1$
denote this composition.

The boundary $\bdy(M_0 \# \s\cpb)$ is a $\S_h$-bundle
over $S^1$ with monodromy
$H=D_{\g_1}\cdots D_{\g_\mu}$,
which is isotopic to the identity.
In fact, the monodromy for the bundle on
$\bdy(M_0 \# \s\cpb)$ has
an additional $\s$
Dehn twists about curves which bound disks, which are of course
isotopic to the identity and can therefore be ignored.
The image of $\omega$ is an $S^2$-bundle over $S^1$
whose monodromy is induced from the word
$D_{\g_1}\cdots D_{\g_\mu}$.
Restricted to each fiber of $\bdy(M_0 \# \s\cpb)$,
$\omega$ agrees with $\pi.$
(As in the strictly nonseparating case of
Section 3, we can see the branching explicitly
in the boundary.
The branched covering $\omega$ now projects each Dehn twist
$D_{\g_i}$ to a disk twist about the arc $\d_i$ for each
nonseparating $\g_i$, and to $D_{\d_i}^2$ for each
separating $\g_i$.
The branch set of $\omega$
appears as the link $L$ in $S^2\times S^1$
that is the
boundary of the
ribbon manifold branch set corresponding to the nonseparating
vanishing cycles,  with an additional
pair of full right-handed twists corresponding to each
separating vanishing cycle. These twists arise from blowing down
$\bdy(S^2\times D^2\# 2\s\cpb)$ to $S^2\times S^1$,
as described in Section 3. Thus, as in Section 3, $L$
is isotopic to a closed braid, and this braid records the motion
of the $2h+2$ branch points of $\pi$
induced from the monodromy
$D_{\g_1}\cdots D_{\g_\mu}$.)

We again have a $\S_h$-bundle equivalence
$\Psi:\bdy M_0 \to \S_h\times S^1$
obtained by expressing $\bdy M_0$ as 
$\S_h\times I/(x,1)\sim (H(x),0)$,
and using the isotopy between $H$ and the
identity to adjust the gluing.
Thus, letting $H_t$ denote this
isotopy, by \cite{bh} we may assume that
$\pi H_t$ projects this to an
isotopy between $\pi H$ and the identity on $S^2$.
This gives a bundle equivalence
$\psi:S^2\times S^1 \to S^2\times S^1$
for which the diagram
$$\begin{CD}
\bdy(M_0\# \cpb)  
@>\Psi>> \S_h\times S^1\\
@V{\omega}VV        @VV{\pi\times id}V\\
S^2\times S^1 @>\psi>> S^2\times S^1
\end{CD}$$
commutes.
Therefore,
$$\vp\cup (\pi\times id): (M_0 \# \s\cpb) \cup_{\Psi^{-1}} (\S_h\times D^2)
\to (S^2\times D^2\# 2\s\cpb) \cup_{\zeta^{-1}} (S^2\times D^2)$$
describes $M\# \cpb$ as a 2-fold cover of 
($S^2$-bundle over $S^2) \# 2\s\cpb$, where $\zeta$ is
the composition of the diffeomorphism
$\bdy(S^2\times D^2\# 2\s\cpb)
\to S^2\times S^1$ and $\psi.$
The branch set is the closed surface obtained by attaching 
the disks
$B\times D^2 \subset S^2\times D^2$ to $L$ via $\zeta^{-1}$.
The proof is completed by noting that 
($S^2$-bundle over $S^2) \# 2\s\cpb \cong
\cp \# (2\s+1) \cpb$.
\end{pf}

\section{An Example}

We illustrate Theorems~\ref{nonsep} and ~\ref{sep} by describing the
branch set in the following example.

In \cite{ma}, Matsumoto showed that $T^2\times S^2 \# 4\cpb$
admits a genus 2 hyperelliptic Lefschetz fibration over $S^2$ with eight
singular fibers. The global monodromy of this fibration
is described by the ordered collection
$(\b_1,\g_0,\b_2,\b_3)^2$
of vanishing cycles,
where $\b_1,\b_2,\b_3,\g_0$
are the curves on $\S_2$ shown in Figure~\ref{matcurves}.
\begin{figure}[h]
$$\epsfbox{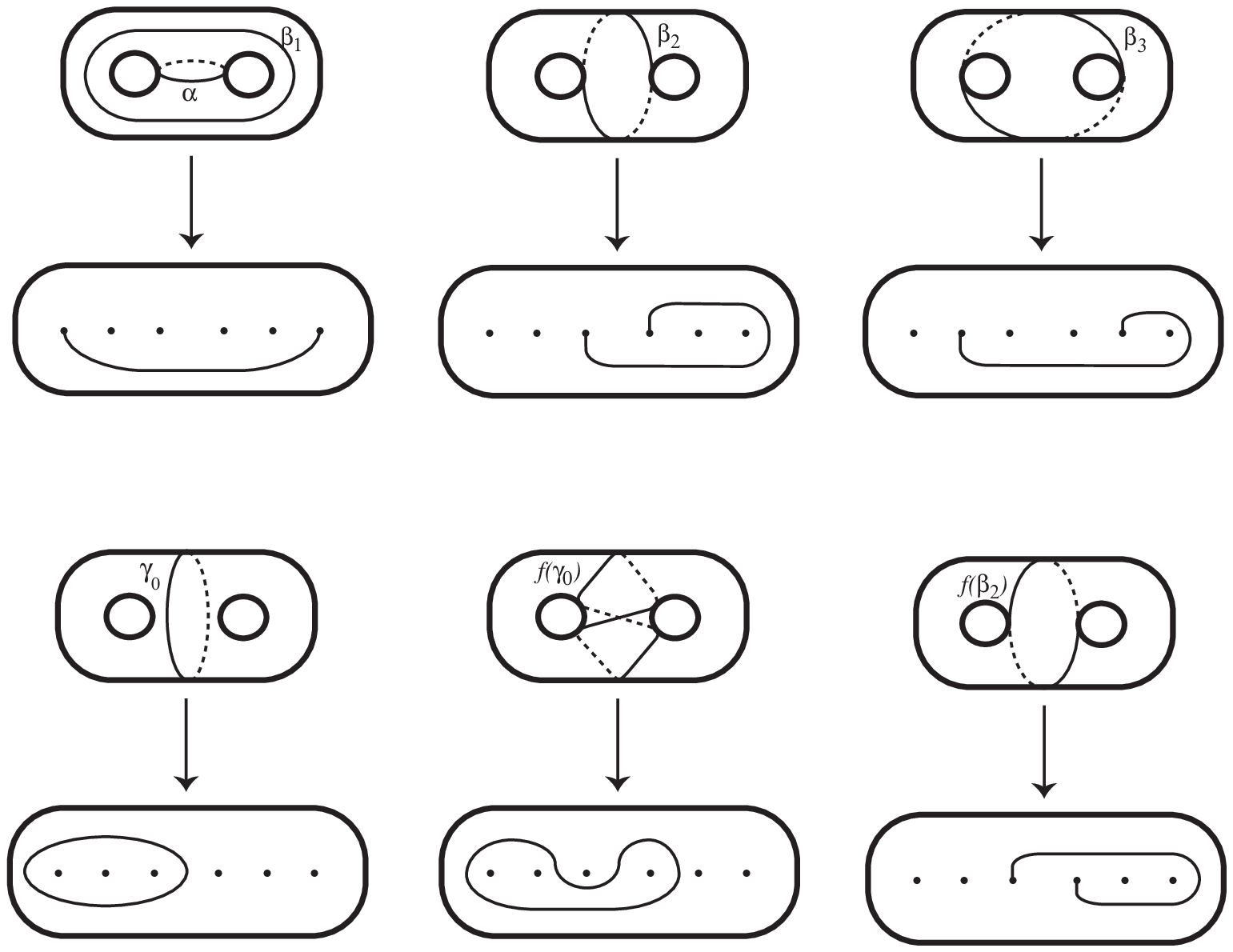}$$
\caption{}
\label{matcurves}
\end{figure}
Let $f=D_{\a},$ where $\a$ is as shown in Figure~\ref{matcurves}.
We define $M$ to be the genus 2 Lefschetz fibration over $S^2$ with
global monodromy
$$(\b_1,\g_0,\b_2,\b_3)^2
(f(\b_1),f(\g_0),f(\b_2),f(\b_3))^2.$$
(This example was purposely constructed to include a
separating curve other than $\g_0$.)

It is clear that $f$ fixes $\b_1$ and $\b_3$.
The curves $\b_1, \b_2, f(\b_2)$ and $\b_3$ are nonseparating, and are 
lifts under $\pi:\S_2\to S^2$  of the indicated arcs
in $S^2.$ Thus, each time they occur as a vanishing cycle,
they contribute a half-twisted 
band (relative to the product)
to the ribbon manifold branch set,
attached according to the corresponding arc.
The curves $\g_0$ and $f(\g_0)$ are separating, 
and are lifts of the indicated simple closed curves in $S^2.$
Whenever they appear as a vanishing
cycle, we attach a $-1$-framed 2-handle to this curve, together
with a $-2$-framed 2-handle attached to its meridian.
Thus we can draw the portion of the branch set residing in
in $S^2\times D^2 \# 8\cpb$ as in Figure~\ref{matbase}.
\begin{figure}[h]
\epsfxsize=5in 
$$\epsfbox{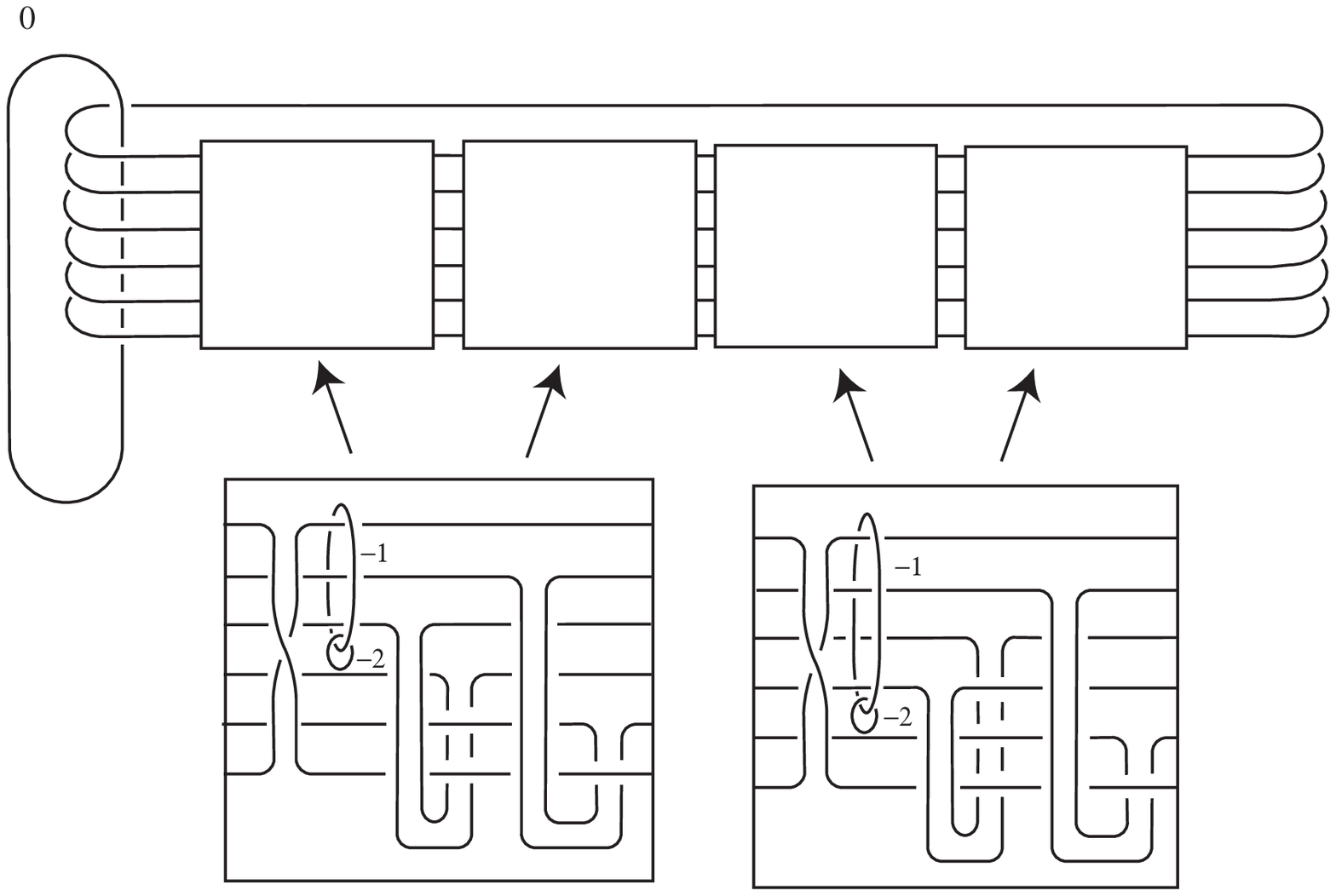}$$
\caption{}
\label{matbase}
\end{figure}
When we close off $S^2\times D^2 \# 8\cpb$
by attaching $S^2\times D^2$, this branch surface
is closed off by six disks attached to the boundary of
the ribbon manifold.
Forming the 2-fold branched cover gives a non-relatively 
minimal Lefschetz fibration on 
$M \# 4\cpb$,
which has $M$
as its relative minimalization.

\section{Deformations versus Resolutions of Branch Surfaces}

In this section, we show that the 
branch surfaces for the branched covers
constructed in this paper are related to
basic desingularization processes in complex algebraic geometry.

Let $f \in {\IC}[z,w]$ be a complex polynomial, and assume that
$f$ has an isolated singularity at the origin. Let
$$V=\{(z,w)\in {\IC}^2 \colon f(z,w)=0\} \cap B_{\varepsilon}^4,$$
where $B_{\varepsilon}^4$ is a 4-ball centered at the origin of 
radius $\varepsilon$, with $\varepsilon$ chosen small enough 
so that $B_{\varepsilon}^4$
contains only the singularity at $\bf{0}$ and so that $V$
intersects $\bdy(B_{\varepsilon}^4)$ transversely.
We view $V$ as a local model for a 
singular portion of a branch surface inside of a 4-manifold,
which we would like to desingularize in order to form a smooth
branched cover.  There are then two standard methods from complex
geometry to do this.  On one hand, we may deform $V$ into
$$V_d=\{(z,w)\in {\IC}^2 \colon f(z,w)=t\} \cap B_{\varepsilon}^4,$$
which will be smooth for $0<|t| \ll \varepsilon.$
On another hand, we may form a resolution $V_r$ of $V$,
which is obtained by repeatedly blowing up $V$ until it becomes
smooth. While any $V_r$ so produced may not be unique, the 
minimal resolution of the surface singularity obtained by forming
the double cover branched over $V_r$ and blowing down any leftover
exceptional curves is unique \cite{bpv}.
See \cite{gs} and \cite{hkk} for further expositions, from a topological
perspective, of these processes.

We demonstrate these procedures with a key example.

{\bf Example.} Let $f(z,w)=z^n+w^{2n},$ where
$n=2g+1$ is odd.
Then $V$ is homeomorphic to a cone on an $(n,2n)$ torus
link. 
The singular point of $V$ is known in algebraic geometry
as an {\it infinitely close n-tuple point.} 
We wish to compare $V_d$ and $V_r,$
as well as the 2-fold covers formed with each
as a branch set.

The deformation $V_d$ of $V$ in $B^4$ is known to be
isotopic to the fibered Seifert surface
of minimal genus in $S^3$ of an $(n,2n)$ torus link \cite{mi}.
See Figure~\ref{def}.
\begin{figure}[h]
$$\epsfbox{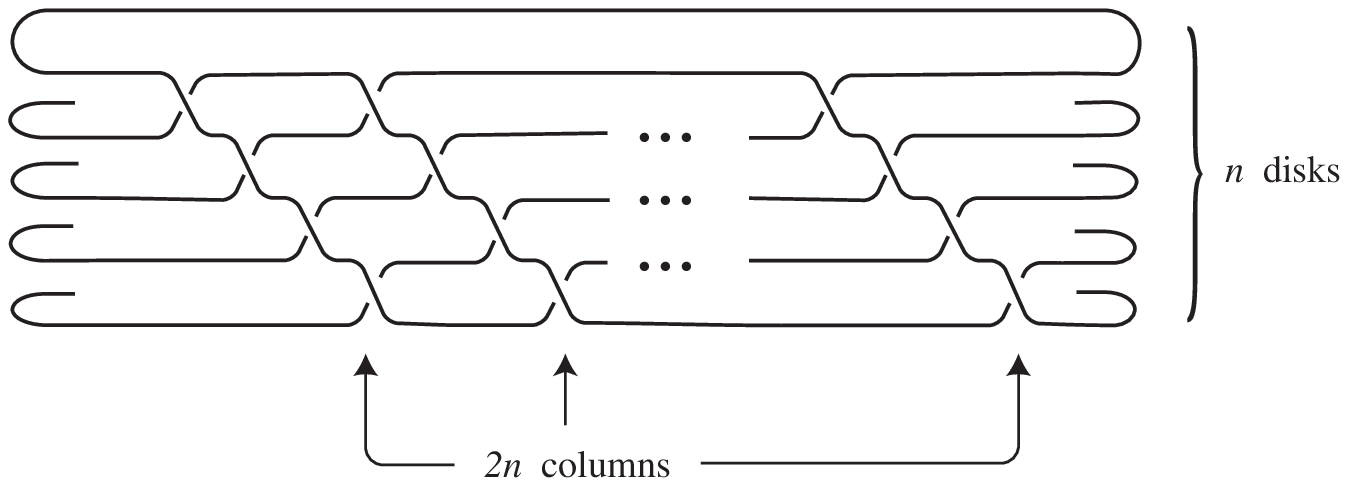}$$
\caption{}
\label{def}
\end{figure}
Figure~\ref{2def} shows a Kirby calculus picture of the
2-fold cover of $B^4$, branched over this surface
(with its interior pushed into $B^4$) \cite{ak}.
\begin{figure}[h]
\epsfxsize=5in 
$$\epsfbox{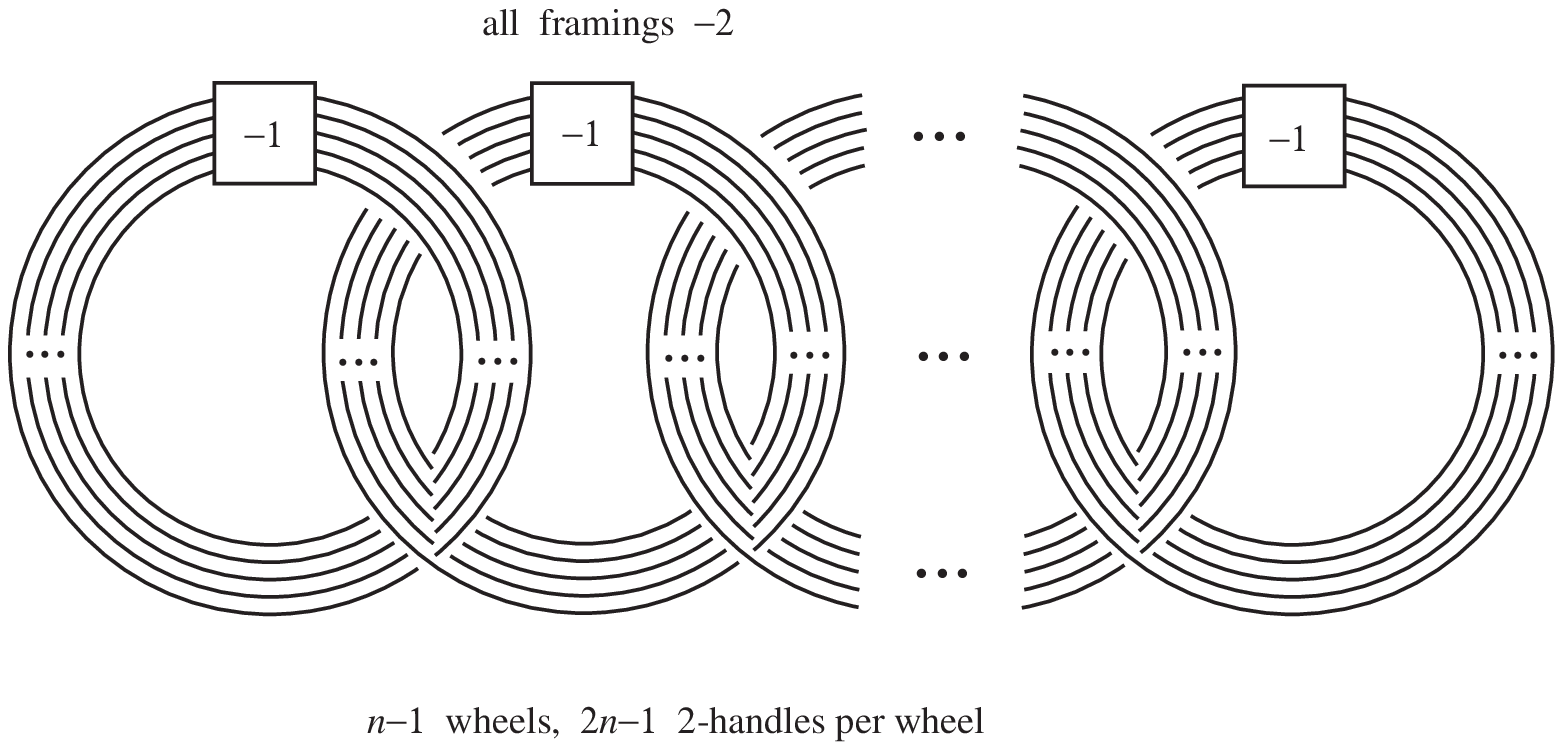}$$
\caption{} 
\label{2def}
\end{figure}
This figure describes a regular neighborhood of the shown
configuration of $(n-1)(2n-1)$ 2-spheres of square $-2$.

We pictorially describe the resolution $V_r$ in 
Figure~\ref{res}.
\begin{figure}[h]
$$\epsfbox{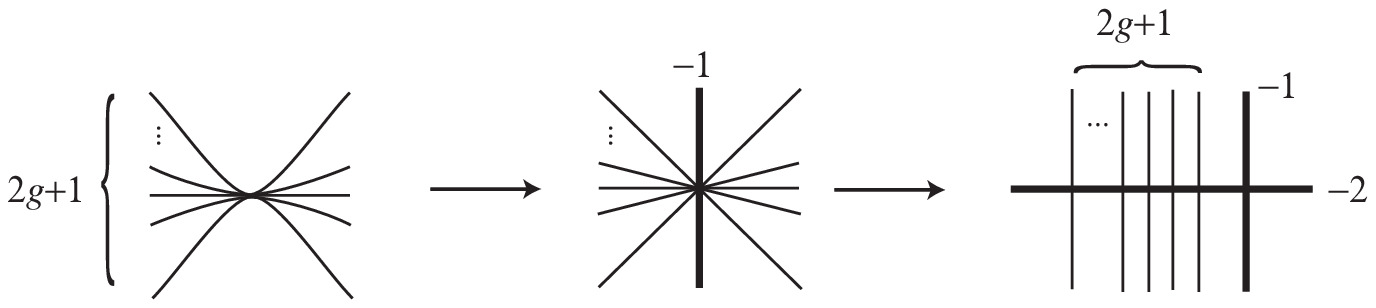}$$
\caption{} 
\label{res}
\end{figure}
(We leave the algebraic computation of the effect of blowing
up in local coordinates, such as in \cite{gs}, to the reader.)
The first blow up shown in Figure~\ref{res} alters 
the singular point of $V$ into
$n$ complex lines meeting transversely in a single point; the 
second blow up then makes the disjoint. 
Since the union of these lines does not represent an even (i.e.
divisible by 2) homology class in
$H_2(B^4 \# 2\cpb),$ we cannot form a double cover
with it as branch set. Instead, we must recall the fact
from the theory of resolutions of complex surface singularities
that when forming 2-fold covers branched over resolutions
of plane curve singularities, we must include in our branch set
any introduced curves of even square \cite{hkk}. 
Let $S$ be 
the embedded 2-sphere of square $-2$ shown in Figure~\ref{res}.
We may then form the 2-fold cover of $B^4\# 2\cpb$
branched over $V_r \cup S.$
The sphere of square $-1$ in Figure~\ref{res} lifts to this cover
as a genus
$g$ surface of square $-2$. This surface intersects a sphere of
square $-1$ obtained as the lift of the sphere of square $-2$
in Figure~\ref{res}. Hence we depict the double cover with
the graph in Figure~\ref{2res}. 
\begin{figure}[h]
$$\epsfbox{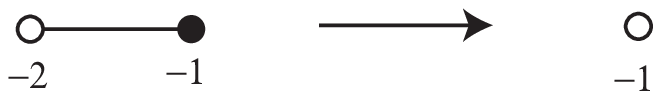}$$
\caption{} 
\label{2res}
\end{figure}
Blowing down the sphere
of square $-1$ gives a regular neighborhood of a genus $g$
surface of square $-1$.

{\Rem As this example demonstrates, double
covers branched over $V_d$ versus $V_r$
may not coincide. It is a 
theorem of Brieskorn that they agree for
complex surface singularities which are
rational double points. A topological proof
of Brieskorn's theorem may be found in \cite{hkk}.}

{\bf A Model Fibration Revisited.}
Let $\a_1, \ldots, \a_{2g}$ and $\g_0$
denote the indicated curves on $\S_h$ in 
Figure~\ref{genush}.
\begin{figure}[h]
\epsfxsize=5in 
$$\epsfbox{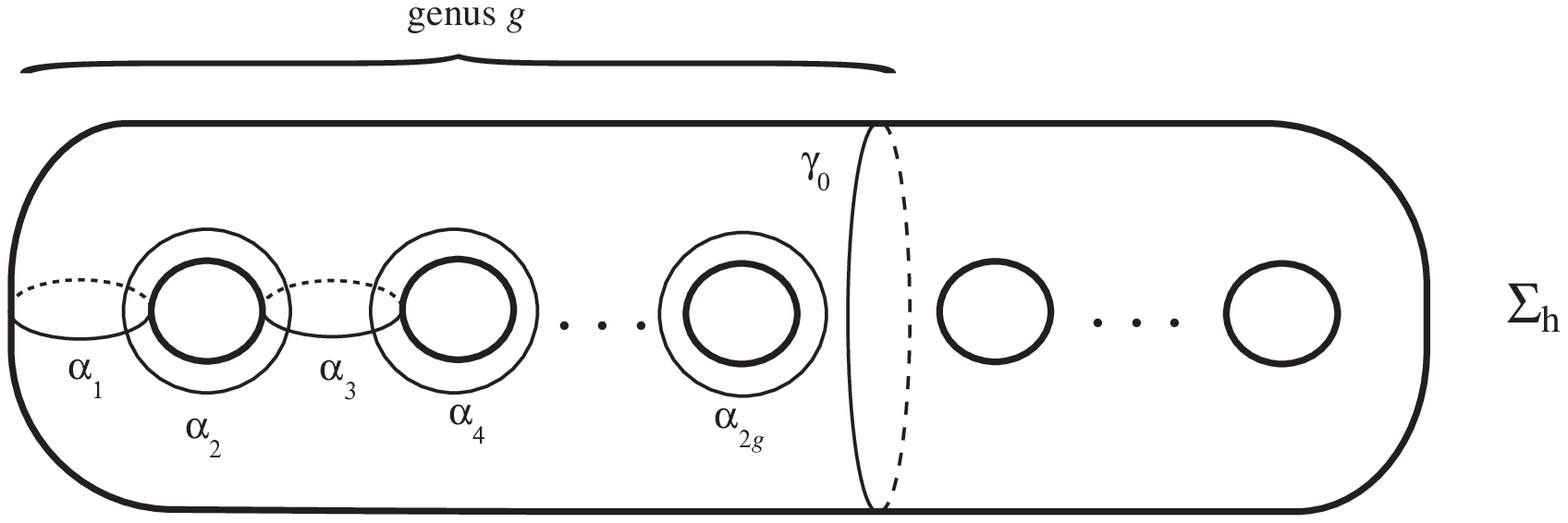}$$
\caption{} 
\label{genush}
\end{figure}
Let $N_0\to D^2$ denote the genus $h$ Lefschetz fibration
of Section 4 with one singular fiber whose vanishing 
cycle is $\g_0$.
Let $P\to D^2$ denote the Lefschetz fibration with global
monodromy given by
$(\a_1,\cdots, \a_{2g})^{2(2g+1)}.$
Since $$(D_{\a_1}\cdots D_{\a_{2g}})^{2(2g+1)}=D_{\g_0}$$ in
${\cal{M}}_h$, $N_0$ and $P$ have the same boundary.
Our goal is to compare the branched coverings that produce
$N_0$ and $P$.

Since all of the curves
$\a_1, \ldots, \a_{2g}$
are nonseparating,
$P$ is obtained as a 2-fold cover of
$S^2\times D^2$ by the methods of Section 3.
The ribbon manifold branch set $B$ can be explicitly
drawn in the boundary, as in Figure~\ref{ribbon}.
\begin{figure}[h]
\epsfxsize=5in 
$$\epsfbox{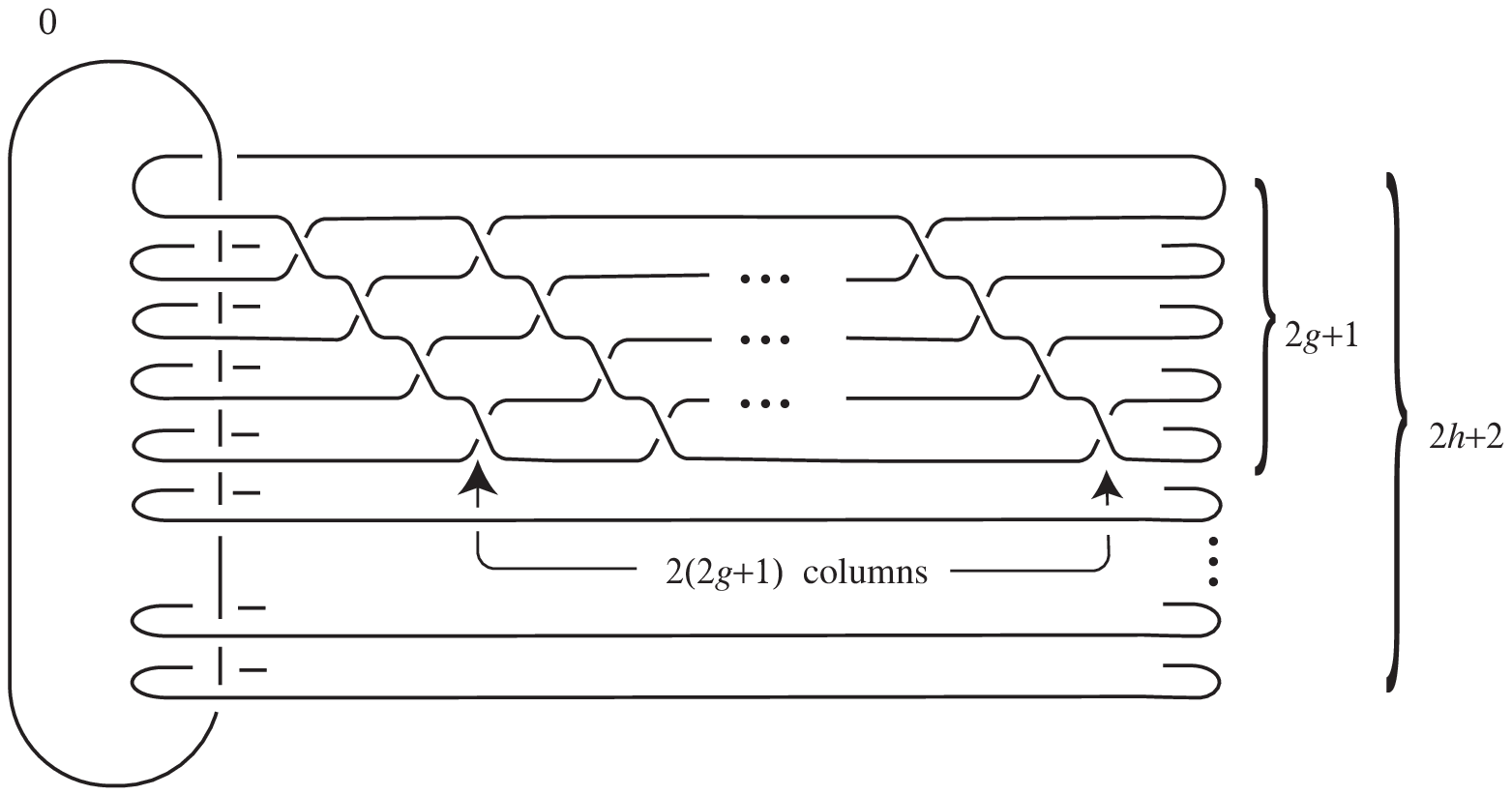}$$
\caption{} 
\label{ribbon}
\end{figure}

Let $D$ denote the 3-ball in
$\bdy(S^2\times D^2)$ shown in Figure~\ref{Bball},
and let $E=D\times I$ denote a 4-ball
collar neighborhood of $D$.
\begin{figure}[h]
\epsfxsize=5in 
$$\epsfbox{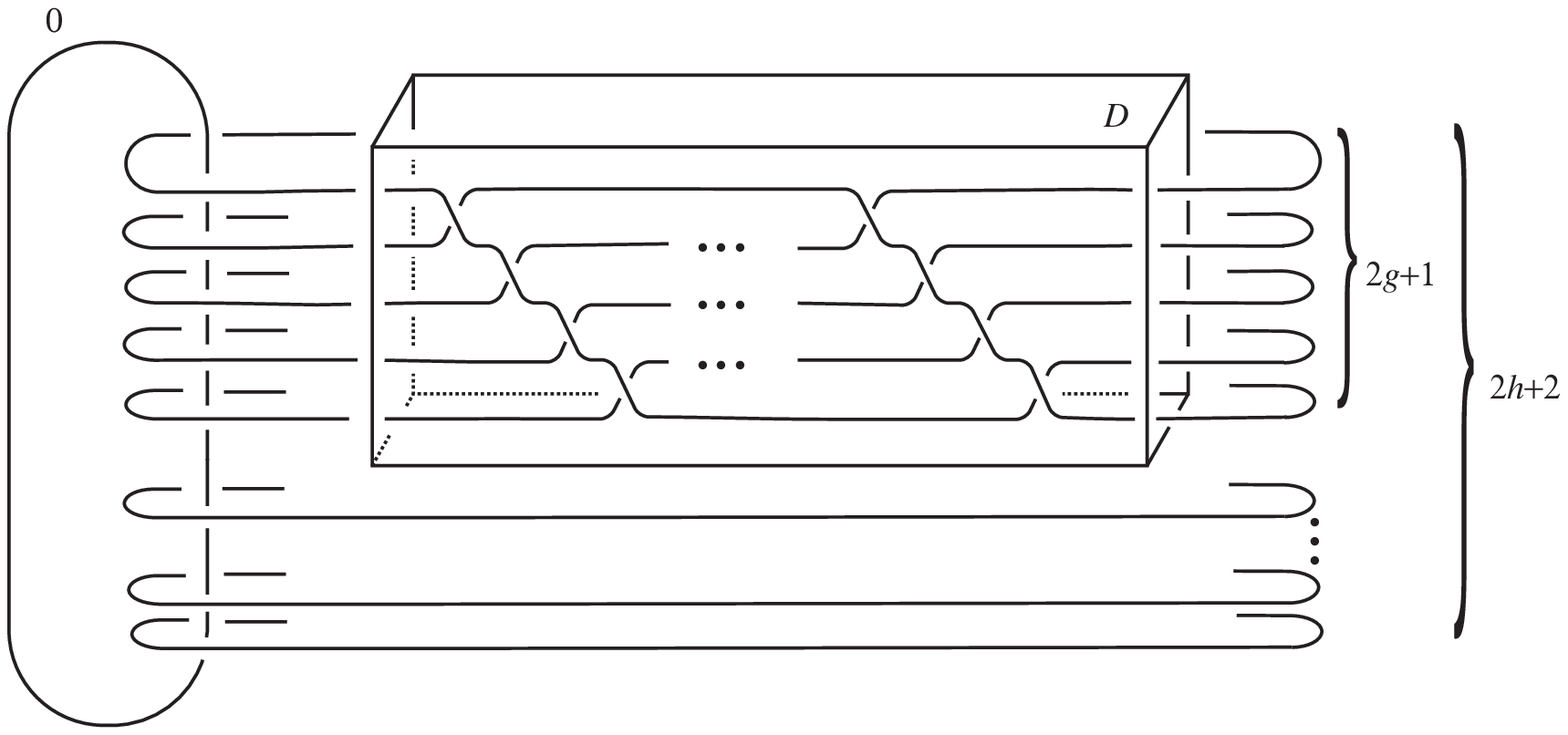}$$
\caption{} 
\label{Bball}
\end{figure}

Then $B\cap E$ is the fibered Seifert surface of minimal genus for a 
$(2g+1,2(2g+1))$ torus link,
so $(E, B\cap E)\cong (B^4,V_d),$
where $V_d$ is the deformation
of an infinitely close $n$-tuple point, with $n=2g+1$.
Hence we may delete this 4-ball, and replace it with 
$(B^4, V)$,
replacing the deformation
$V_d$ with $V$.
This produces an immersed surface $C=(B-V_d)\cup V$.
Alternatively, we may delete $(B^4,V_d),$
and replace it with $(B^4\# 2\cpb, V_r).$
(It is a routine Kirby calculus exercise to verify
that the boundary diffeomorphism
$\bdy(B^4\# 2\cpb) \to S^3$
given by blowing down the final picture in 
Figure~\ref{res} sends $\bdy(V_r)$ to a
$(2g+1,2(2g+1))$ torus link. Thus we may match
up $\bdy(V_r)$ and $\bdy(V_d)$.)
Doing this transforms the manifold in
Figure~\ref{Bball} to the manifold in
Figure~\ref{b2}. By Proposition~\ref{model},
forming the 2-fold branched cover of the latter branched
over $B^{\prime}=(B-V_d)\cup V_r\cup S$ produces a non-relatively
minimal Lefschetz fibration 
whose relative minimalization is
$N_0\to D^2$.
We have shown the following.
{\Prop \label{defvsres}
\begin{itemize}
\item[(a)] 
The branch surface for the branched covering
$N_0\# \cpb \to S^2\times D^2\# 2\cpb$
is the resolution of the immersed surface $C \cup S$ in $S^2\times D^2$.
\item[(b)] The branch surface for the branched covering
$P \to S^2\times D^2$ is the deformation of the immersed 
surface $C$ in $S^2\times D^2$.
\end{itemize}}

The proof of Theorem~\ref{sep} used the action
of ${\cal H}_h$ on separating curves to transfer the model
of a branched cover giving $\g_0$
as a vanishing cycle to give any symmetric separating curve
as a vanishing cycle. Similar arguments give the following.

{\Cor \label{allres}
Let $M\to S^2$ be a relatively minimal genus $h$ Lefschetz fibration
whose monodromy includes $\s$ separating vanishing cycles.
Then 
the branch surface for the branched covering
$M \# \s\cpb \to \cp \# (2\s+1) \cpb$
is the resolution of an immersed surface
in $S^2\times S^2$ or \snx with only infinitely close
n-tuple point singularities.}

\begin{pf}
If $\g_i$ is any separating vanishing cycle which
bounds a genus $g$ surface, then we have as in the 
proof of Theorem~\ref{sep} a commutative diagram
$$ \begin{CD}
(\S_h\times D^2, \g_0) @>\tilde{f}\times id>> (\S_h\times D^2,\g_i)\\
@V{\pi\times id}VV          @VV{\pi\times id}V\\
(S^2\times D^2,\d_0)  @>f\times id>> (S^2\times D^2,\d_i)
\end{CD}$$
which extends to a commutative diagram
$$\begin{CD}
N_0\# \cpb @>\widetilde{F}>> N_i\# \cpb\\
@V{\vp_0}VV          @VV{\vp_i}V\\
S^2\times D^2\# 2\cpb  @>F>> S^2\times D^2\# 2\cpb.
\end{CD}$$
The portion of the branch set constructed in the proof of
Theorem~\ref{sep} to produce $\g_i$ as a separating 
vanishing cycle can be described as the result of 
replacing $(f\times id)(B^4, V)$ with
$F(B^4\# 2\cpb, V_r)$.
\end{pf}

{\Rem
In general, a right-handed Dehn twist about
a separating curve on $\S_h$ can always be expressed as a
product of right-handed Dehn twists about nonseparating
curves. Thus, we can always consider the operation on Lefschetz
fibrations of removing the neighborhood of a singular fiber
with a separating vanishing cycle, and replacing it with 
a fibration given by the appropriate collection of nonseparating
vanishing cycles (and vice versa). 

We can isolate this operation to smaller neighborhoods,
as follows.
When these fibrations
are hyperelliptic and viewed as branched covers, 
we have shown that this operation
can be phrased as replacing a submanifold of the 
base diffeomorphic to $(B^4 \# 2\cpb, V_r)$
with $(B^4, V_d)$, and forming 2-fold
covers.
Lifting these submanifolds to the branched cover, 
these operations involve trading back and forth 
(after blowing down exceptional
curves)
the regular 
neighborhood of
an embedded genus $g$ surface of square $-1$ 
as in Figure~\ref{2res} with the regular neighborhood of 
a configuration of
embedded 2-spheres as in Figure~\ref{2def}.}


\begin{thebibliography}{HKK}

\bibitem[AK]{ak}
S. Akbulut and R. Kirby, 
{\em Branched covers of surfaces in 4-manifolds,}
Math. Ann. {\bf 252} (1980), 111-131.

\bibitem[BPV]{bpv}
W. Barth, C. Peters, A. Van de Ven,
{\em Compact complex surfaces,} 
Ergebnisse der Mathematik, Springer-Verlag, Berlin, 1984.

\bibitem[B]{b} 
J. Birman, {\em Braids, links, and mapping class groups,}
Ann. of Math.Stud. {\bf 82}, Princeton Univ.
Press, Princeton, NJ, 1975.

\bibitem[BH]{bh} J. Birman, H. Hilden,
{\em On the mapping class groups of closed surfaces
as covering spaces,}
Advances in the theory of Riemann surfaces,
Ann. of Math.Stud. {\bf 66}, Princeton Univ.
Press, Princeton, NJ (1971), 81-115.

\bibitem[D]{d} S. Donaldson, {\em Lefschetz fibrations in symplectic
geometry,}
Doc. Math. J. DMV Extra Volume ICM
(1998) II, 309-314.

\bibitem[GS]{gs} 
R. Gompf and A. Stipsicz,
{\em An introduction to 4-manifolds and Kirby calculus,}
book in preparation.

\bibitem[HKK]{hkk}
J. Harer, A. Kas, and R. Kirby,
{\em Handlebody decompositions of complex surfaces,}
Mem. Amer. Math. Soc. {\bf 62} (1986), number 350.

\bibitem[K]{k} 
A. Kas, {\em On the handlebody decomposition associated
to a Lefschetz fibration,} Pacific J. Math.
{\bf 89} (1980), 89-104.

\bibitem[Ma]{ma}
Y. Matsumoto, {\em Lefschetz fibrations of genus two - a topological
approach,} Proceedings of the 37th Taniguchi Symposium on
Topology and Teichm\"{u}ller Spaces, ed. Sadayoshi Kojima
et al., World Scientific (1996) 123-148.

\bibitem[Mi]{mi} J. Milnor, 
{\em Singular points of complex hypersurfaces,}
Ann. of Math. Stud. {\bf 61},
Princeton Univ. Press, Princeton, NJ, 1968.

\bibitem[Mo]{mo} J. M. Montesinos, {\em 4-manifolds,
3-fold coverings, and ribbons,} Trans. Amer. Math. Soc.
{\bf 245} (1978), 453-467.

\bibitem[Ms]{ms} 
B. Moishezon,
{\em Complex surfaces and connected sums of complex
projective planes,} Lecture Notes in Math.
{\bf 603}, Springer, New York, 1977.

\bibitem[P1]{p1} 
U. Persson,
{\em Chern invariants of surfaces of general
type}, Compositio Math. {\bf 43} (1981), 3-58. 

\bibitem[P2]{p2} 
U. Persson,
{\em Genus 2 fibrations revisited,}
Lecture Notes in Math. {\bf 1137},
Springer, New York (1992), 133-144. 

\bibitem[ST]{st} B. Siebert and G. Tian,
{\em On hyperelliptic $C^{\infty}$-Lefschetz fibrations
of four-manifolds,} preprint
{\tt math.GT/9903006}.

\bibitem[S]{s} I. Smith, {\em Symplectic geometry of
Lefschetz fibrations,}
Dissertation, Oxford University, 1998.


\end{thebibliography}
\end{document}